\documentclass{article}

\usepackage[T1]{fontenc}
\usepackage{booktabs}
\usepackage{graphicx}
\usepackage{multirow}
\usepackage{url}
\usepackage[cmex10]{amsmath}
\usepackage{amssymb}
\usepackage{natbib}
\usepackage{authblk}
\usepackage{hyperref}
\usepackage{algorithmic}
\usepackage{bm}
\usepackage{amsthm}
\allowdisplaybreaks

\newtheorem{theorem}{Theorem}[section]
\newtheorem{remark}{Remark}

\newtheorem{corollary}{Corollary}[section]
\newtheorem{lemma}{Lemma}[section]

\usepackage{geometry}
\begin{document}
\title{\Large{\bf Law of the Iterated Logarithm and Model Selection Consistency for GLMs with Independent and Dependent Responses}}

\author{{Xiaowei Yang}$^1$,~~~~{Shuang Song}$^{2,3}$,~~~~{Huiming Zhang}$^{^4}$\footnote{Correspondence author. Email: \texttt{zhanghuiming@pku.edu.cn} (Huiming Zhang)}\\
~~\\
{\small{1.School of Mathematics and Statistics, Chaohu University, Hefei, Anhui, China \\
2.Center for Statistical Science, Tsinghua University, Beijing, China\\
3.Department of Industrial Engineering, Tsinghua University, Beijing, China\\
4.School of Mathematical Sciences, Peking University, Beijing, China\\
}}}
\date{}  
\maketitle

\begin{abstract}
We study the law of the iterated logarithm (LIL) for the maximum likelihood estimation
of the parameters (as a convex optimization problem) in the generalized linear models with independent or weakly dependent ($\rho$-mixing, $m$-dependent) responses under mild conditions. The LIL is useful to derive the asymptotic bounds for the discrepancy between the empirical process of the log-likelihood function and the true log-likelihood. As the application of the LIL, the strong consistency of some penalized likelihood based model selection criteria can be shown. Under some regularity conditions, the model selection criterion will be helpful to select the simplest correct model almost surely when the penalty term increases with model dimension and the penalty term has an order higher than $O({\rm{loglog}}n)$ but lower than $O(n)$. Simulation studies are implemented to verify the selection consistency of BIC.

\end{abstract}

\maketitle

\section{INTRODUCTION}

Originating from \cite{Nelder1972}, generalized linear models (GLMs) are remarkable
and synthetic extensions of linear models. GLMs are often classified into two classes in references. The first type is the GLM with a natural link function (canonical link function), such as the binomial regression (logistic regression) and the Poisson regression. GLMs endowed with non-natural link functions become the second type of GLMs, including the probit model and the negative binomial regression, which are more complex to analyze, see \cite{Fahrmeir85}, \cite{Chen11}. Recently, GLMs have become the very popular regression models in the big data era, see \cite{Efron2016}.

Under some regularity conditions, asymptotic normality for maximum likelihood estimator (MLE) in GLMs with both natural and non-natural link functions, was established by \cite{Fahrmeir85}. However, asymptotic normality is a type of convergence which is weaker than the strong limit theorems. A lot of efforts have been done in studying strong consistency in terms of the law of iterated logarithm (LIL). In linear models, \cite{Lai1982} obtained a general LIL for weighted sums of the independent regression noises with zero means and common variance, which can be applied to least squares estimation. \cite{He1995} proved the LIL for a general class of M-estimators and then derived the Bahadur representation. \cite{Fang1998} studied the LIL for the MLE of general nonhomogeneous Poisson processes. If there are measurement errors in covariates, under some regularity conditions, \cite{Miao2011} established the LIL for the least squares estimation in simple linear models.


For many practical regression problems, redundant or irrelevant covariates will come into the models. If those irrelevant covariates are enclosed, the efficiency of estimators will more or less be impaired. It is indispensable to do variable selection after obtaining the regression coefficients in scientific analysis. Usually, when we consider proving the consistency of model selection like the BIC, weak or even strong consistency of the estimated coefficient is not enough if convergence rate cannot be obtained. Thus, it is worth concerning the rate of strong consistency of the estimator as is the LIL. Under i.i.d. noise assumptions (may be non-normal), \cite{Rao1989} applied the LIL to the selection consistency of linear models with non-normal noise. Similarly,  \cite{Wu1999} studied the Huber's M-estimator for linear models.

This work focuses on studying a set of the more flexible GLMs where more link functions and exponential family responses are permitted. In some simple GLMs, canonical links do not always provide the best fit. For example, as the non-canonical link GLM, the negative binomial regression can model over-dispersed count data while the Poisson regression (a GLM with the canonical link function) requires equal dispersion, which is inapplicable in practice. Generally, there is no apriori reason to explain why a canonical link should be used, and in many cases a non-canonical link is more suitable (\cite{McCullagh1989} and \cite{Czado00}). Note that \cite{Qian2006}(merely binomial regression) only considered some special cases of GLMs. Our results are the extensions of \cite{Qian2006} but with slightly differences (see Remark1, (H.5) and (H.6)). Meanwhile, we require fewer assumptions on the Fisher information matrix of true regression parameters. 

Let $\beta_{0}:=(\beta_{01},\ldots,\beta_{0p})^{t}$ be the true value of regression parameter $\beta$ and $\|\cdot\|$ be the Euclidian norm. In order to get the model selection consistency, we need to show the LIL for the MLE $\hat{\beta}$ of GLMs. Assume that there exists a constant $d>0$ such that
\begin{equation}\nonumber
\mathop {\limsup }\limits_{n \to \infty } \frac{\parallel\hat{\beta}-\beta_{0}\parallel}{\sqrt{n^{-1}\log\log n}}=d \qquad a.s.,
\end{equation}
under some conditions (see Section 3.3 for details). Then we have the convergence rate
$
\parallel\hat{\beta}-\beta_{0}\parallel=O(\sqrt{n^{-1}\log\log n})~a.s..
$ Hence the strong consistency of $\hat{\beta}$ is directly implied by the LIL.

Moreover, this work gives a general consideration of model selection consistency in GLMs with non-natural link functions and the responses are allowed to be weakly dependent including $\rho$-mixing and $m$-dependent. Our contributions are to extend previous works under some mild conditions, which lead to desired performance of model selection, while the previous studies about GLMs seldom consider the non-natural link function and weakly dependent responses. Our generalization covers a wide range of GLMs.

This paper is organized as follows. In Section 2, a review of exponential family GLMs is presented and the problem of estimating the coefficients in GLMs by weighted scores equation is discussed. In section 3, under regularity assumptions, we give the LIL for the weighted scores estimates for independent or weakly dependent responses. Checking the model selection consistency is equivalent to evaluating whether the order of penalty term is between $O({\rm{loglog}}n)$ and $O(n)$. For independent or weakly dependent exponential family responses, it is easy to show the model selection consistency by applying the LIL proposed in Section 3. In section 4, we run the simulation for 500 times, and exemplifies that the BIC is consistent while the AIC is not consistent. The detailed proofs of the theorems and lemmas are demonstrated in APPENDIX.

\section{GENERALIZED LINEAR MODELS}
We expect readers to be familiar with exponential families as prerequisites. Here we give a brief review that makes our exposition self-contained.

Let $\{(y_i,x_i)\}_{i=1}^{n}$ be $n$ independent sample data pairs. Let ${x_i}^{,}s\in \mathbb{R}^{p}$ be fixed covariates and assume that the responses ${y_i}^{,}s$ follow the density of the exponential families
\begin{equation}\label{eq-glm}
dF(y_i)=c(y_i)\exp\{y_i\theta_i-b(\theta_i)\}d\mu(y_i),c(y_i)>0,i=1,\ldots,n
\end{equation}
where ${\theta _i} \in \Theta  = \{\eta :\int c (y)\exp \{ y\eta \} \mu (dy) < \infty \}$, $c(y_i)$ is free of $\theta_i$ and $\mu(\cdot)$ is the dominated measure.

Let $\dot{b}(\cdot)$ and $\ddot{b}(\cdot)$ be the first and second derivatives of $b(\cdot)$, respectively. A notable property of the expectation and variance of ${y_i}^{,}s$ in \eqref{eq-glm} is $E(y_i)=\dot{b}(\theta_i)$ and $Var(y_i)=\ddot{b}(\theta_i)$. Note that the first and second derivatives of $b(\cdot)$ are bounded in any compact set $K \subset \Theta$.  For more details of exponential families, please see \cite{Brown1986}.

For GLMs, we assume that the effect of covariates $x_i$ on responses $y_i$ can be observed by the corresponding regression coefficient $\beta$, and the relationship between $\theta_i$ and $x_{i}^{t}\beta$ is expressed by the following \emph{relation function}
$$
\theta_i=u(x_{i}^{t}\beta_0).
$$
Therefore, $\theta_{i}$ is determined by the unknown $\beta$ and the given $x_i$. Here we suppose that $x_i$'s are non-random vectors of dimension $p$. The $\beta\in \mathbb{R}^{p}$ is called regression coefficient which will be estimated.

Further, we have
\[E({y_i}|{x_i}): = {\mu _i} = \dot b({\theta _i}) = \dot b(u(x_i^t\beta_0 )),i = 1,2, \cdots ,n.\]
Let $g$ be the \emph{link function} such that $g(\mu_i)=x_{i}^{t}\beta_0$. Since $\mu_{i}=g^{-1}(x_{i}^{t}\beta_0)$ with $\mu_{i}=\dot{b}(u(x_{i}^{t}\beta_0))$, we immediately acquire the expression $u(t)=\dot{b}^{-1}(g^{-1}(t))$.

The likelihood function of $\{(y_i,x_i)\}_{i=1}^{n}$ is the product of $n$ terms in (\ref{eq-glm}). Taking the logarithm, we get the log-likelihood function
\begin{equation}\label{eq-log-like}
l_{n}(\beta):=\sum_{i=1}^{n}[y_{i}u(x_{i}^{t}\beta)-b(u(x_{i}^{t}\beta))].
\end{equation}
Then the score function is defined by
\begin{equation}\label{eq-score}
S_{n}(\beta):=\frac{\partial l_{n}(\beta)}{\partial \beta}=\sum_{i=1}^{n}x_{i}\dot{u}(x_{i}^{t}\beta)[y_{i}-\dot{b}(u(x_{i}^{t}\beta))].
\end{equation}

For more discussion of GLMs, readers are suggested to see \cite{McCullagh1989}, \cite{Fahrmeir2001}, \cite{Chen11}.

\subsection{Examples of Non-natural Link}

For GLMs with the \emph{natural link function} $g(t)=\dot{b}^{-1}(t)$, we have the identity function $u(t)=t$. Thus the first term in (\ref{eq-log-like}) is a linear function of $x_{i}^{t}\beta$ and the score function is
$$
S_{n}(\beta)=\sum_{i=1}^{n}x_{i}\dot{u}(x_{i}^{t}\beta)[y_{i}-\dot{b}(u(x_{i}^{t}\beta))]=\sum_{i=1}^{n}x_{i}[y_{i}-E(y_{i})].
$$
Furthermore, we illustrate some GLMs with non-natural link functions.

\textbf{1. Probit Model}

Let $\Phi(t)$ be the cumulative distribution function of $N(0,1)$. Suppose the $n$ independent observed data sets $\{y_i,x_i\}_{i=1}^n$ are from a logit model, and $y_i$'s are independent Bernoulli distributed with $P(y_i=1|x_i;\beta)={\Phi(x_i^t\beta )},~i=1,2, ..., n$. The average log-likelihood function for data $\{y_i,x_i\}_{i=1}^n$ is
\[l(y, \beta) = \frac{1}{n}\sum\limits_{i = 1}^n {\{ {y_i}\log \Phi (x_i^t\beta ) + (1 - {y_i})\log [1 - \Phi (x_i^t\beta )]\} } .\]

Like the logistic regression, the $u( \cdot )$ and $b( \cdot )$ in the probit model follow parametric equations:
$$
u(x_{i}^{t}\beta)=\log[\Phi(x_{i}^{t}\beta)/(1-\Phi(x_{i}^{t}\beta))],~~b(u(x_{i}^{t}\beta))=\log(1-\Phi(x_{i}^{t}\beta)),
$$
where $\Phi(t)$ is the cumulative distribution function of the standard normal distribution.

\textbf{2. Negative Binomial Regression}

It is  common that negative-binomial GLMs are more plausible than Poisson regression for modelling overdispersion count data. Negative binomial regression (NBR) assumes that the overdispersed response data are modelled by two-parameter distribution:
\begin{eqnarray}\nonumber
f(y_{i}|\theta,\mu_{i})=\frac{\Gamma(\theta+y_{i})}{\Gamma(\theta)y_{i}!}(\frac{\mu_{i}}{\theta+\mu_{i}})^{y_{i}}(\frac{\theta}{\theta+\mu_{i}})^{\theta}, \quad i=1,2,\cdots,n.
\end{eqnarray}
where $\theta$ is known (can be estimated previously), $E(y_{i}|x_{i})=\mu_{i}$, $Var(y_{i}|x_{i})=\mu_{i}+\frac{\mu_{i}^{2}}{\theta}$.

Suppose that the relationship between the mean parameter and covariates is given by $\log(\mu_{i})=x_{i}^{t}\beta$. Then the logarithm of the maximum likelihood function for NBR is
\begin{eqnarray}\nonumber
l(y;\beta)&=&\log\{\prod\limits_{i=1}^{n}f(y_{i}|\theta,\mu_{i})\}\propto\sum\limits_{i=1}^{n}\{y_{i}\log\mu_{i}-(\theta+y_{i})\log(\theta+\mu_{i})\}\nonumber\\
&=&\sum\limits_{i=1}^{n}[y_{i}x_{i}^{t}\beta-(\theta+y_{i})\log(\theta+e^{x_{i}^{t}\beta})].\nonumber
\end{eqnarray}
Thus the connection of $u( \cdot )$ and $b( \cdot )$ for NBR is
$$
u(x_{i}^{t}\beta)=x_{i}^{t}\beta-\log(\theta+\exp\{x_{i}^{t}\beta\}),~b(u(x_{i}^{t}\beta))=\theta\log(\theta+\exp\{x_{i}^{t}\beta\})
$$
where $\mu_{i}=\exp\{x_{i}^{t}\beta\}$.

NBR plays an important role in modern applications relating to count data, see \cite{Zhang2017} and references therein.

\subsection{Weighted scores equations}

When estimating regression coefficients, a well-known robust approach is the weighted-likelihood method which assigns a sequence of weights to perturb the contribution of each sample in the log-likelihood function. Based on some weight functions (see \cite{Markatou98}), the simultaneous changes in the weights of $n$ samples producing the robust estimators for GLMs. It is a method of mitigating the influence of outliers (or say the leverage points, see Pages 59 of \cite{Vaart1998}).

If the log-likelihood (\ref{eq-log-like}) is replaced by the following weighted log-likelihood
\begin{eqnarray}\label{wlog}\nonumber
{l_n (\beta )}:={\sum\limits_{i = 1}^n {w_i [y_i u(x_i^t \beta ) - b(u(x_i^t \beta ))]} },
\end{eqnarray}
the maximum likelihood estimation for the weighted-likelihood is
\begin{equation}\label{eq:mlog}
\hat{\beta}_{n}=\arg\max_{\beta\in \mathbb{R}^{p}}\{\sum_{i=1}^{n}w_{i}[y_{i}u(x_{i}^{t}\beta)-b(u(x_{i}^{t}\beta))]\}.
\end{equation}
Here $\{w_{i}\}_{i=1}^{n}$ are any uniformly bounded positive weights that are independent of response $y_i$ for all $i$ (The weights may depend on $\{x_{i}\}_{i=1}^{n}$).

Then the weighted MLE is a robust estimation for some suitable weights $\{w_{i}\}_{i=1}^{n}$. If $w_{i}\equiv$~constant for all $i$, then the unweighted MLE is exactly the common GLMs, which may not have the robust property.

Consider vector derivative with respect to $\beta$ by letting $\frac{\partial l_{n}(\beta)}{\partial\beta}=\{\frac{\partial l_{n}(\beta)}{\partial\beta_{1}},\cdots,\frac{\partial l_{n}(\beta)}{\partial\beta_{p}}\}^{t}$. Then weighted score is given by
\begin{eqnarray}\label{log-l5}\nonumber
S_{n}(\beta):=\frac{\partial l_{n}(\beta)}{\partial \beta}&=&\sum\limits_{i=1}^{n}w_{i}x_{i}\dot{u}(x_{i}^{t}\beta)[y_{i}-\dot{b}(u(x_{i}^{t}\beta))].
\end{eqnarray}
In particular, we have $S_{n}(\beta_0)=\sum\limits_{i=1}^{n}w_{i}x_{i}\dot{u}(x_{i}^{t}\beta_0)[y_{i}-E(y_{i})]$.

From $\frac{\partial l_{n}(\beta)}{\partial\beta}$, assume that $\hat{\beta}$ is a unique solution by the \emph{weighted scores equation} $S_{n}(\beta)=0$. Note that the Hessian matrix of $\beta$ is
\begin{eqnarray}\label{Hessian-non}
\frac{\partial^{2}l_{n}(\beta)}{\partial\beta\partial\beta^{t}}=\sum\limits_{i=1}^{n}w_{i}\{x_{i}\ddot{u}(x_{i}^{t}\beta)x_{i}^{t}[y_{i}-\dot{b}(u(x_{i}^{t}\beta))]
-x_{i}\dot{u}^2(x_{i}^{t}\beta)\ddot{b}(u(x_{i}^{t}\beta))x_{i}^{t}\},
\end{eqnarray}
and the Fisher information at $\beta$ is
\begin{eqnarray}\label{Fisher-non}\nonumber
I_{n}(\beta):=-E\frac{\partial^{2} l_{n}(\beta)}{\partial\beta\partial\beta^{t}}=-\sum\limits_{i=1}^{n}w_{i}\{x_{i}\ddot{u}(x_{i}^{t}\beta)x_{i}^{t}[E y_{i}-\dot{b}(u(x_{i}^{t}\beta))]
-x_{i}\dot{u}^2(x_{i}^{t}\beta)\ddot{b}(u(x_{i}^{t}\beta))x_{i}^{t},
\end{eqnarray}
where the last equality is due to $E(y_{i})=\dot{b}(u(x_{i}^{t}\beta))$. In particular, we have $$I_{n}(\beta_0)=\sum\limits_{i=1}^{n}w_{i}\dot{u}^{2}(x_{i}^{t}\beta_0)\ddot{b}(u(x_{i}^{t}\beta_0))x_{i}x_{i}^{t}.$$

It is easy to see that the Fisher information $I_{n}(\beta)$ is semi-positive definite. Then, $\hat{\beta}$ makes the likelihood function get the maximum.

\section{MODEL SELECTION CONSISTENCY}
\label{sec2}

\subsection{Preliminaries and Notations}
In order to examine the predictability of the model, in this section we discuss the model selection methods by including or excluding variables. To be specific, for $i=1, \cdots ,n$, the question is ``Which optimal subset of $x=({x_{i1}}, \cdots ,{x_{ip}})^{t}$ will enter into the regression function $E(y_i|x_i)$ by some model selection approaches?''. In the existing literatures, there are several commonly used model selection criteria such as AIC (Akaike Information Criterion, \cite{Akaike1973}), BIC (Bayesian Information Criterion, \cite{Schwarz1978}) and SCC (Stochastic Complexity Criterion, \cite{Rissanen1989}), whose theoretical background is rooted in the information theory or the Bayesian analysis.

To answer the question, we now build a mathematical framework to get optimal sub-model. Let $\alpha$ be a subset of $\{1,2,\ldots,p\}$ and define $\beta(\alpha)$ (or $x_{\alpha}$) to be the sub-vector of $\beta$ (or $x$) indexed by the integers in $\alpha$. The dimension of $\beta({\alpha})$ is denoted by $p_{\alpha}$. Therefore, for simplicity, we say that the sub-model
corresponding to $\alpha$ is called the $\alpha$ sub-model or the candidate model. It is given by the mean of the predictor
\begin{eqnarray}\label{submodel}\nonumber
{\mu _i}: = E({y_i}|{x_i}) = \dot b(u(x_{i\alpha }^t\beta (\alpha ))),~~i = 1,2, \cdots ,n.
\end{eqnarray}
Note that $\alpha$ sub-model is not necessarily a correct model
in a manner that $E({y_i}|{x_i})$ is not constantly and precisely equal to $\dot b(u(x_{i\alpha }^t\beta (\alpha )))$. Let ${\cal A}$ be the collection of all subsets of ${\{ 1,2,\cdots,p\} }$. Then there are $|{\cal A}|=2^{p}$ candidate models for screening. If $\beta_{0i}=0$ for all i not in alpha then model alpha is called "correct". The set of all correct models is denoted $\Gamma_c$.  The set $\alpha_0 = \{i:\beta_{0i} \neq 0\}$ is the smallest correct model.  The set of all models which are not correct is denoted $\Gamma_w$. Then
model $\alpha$  is called a correct
model, denoted by
$$\Gamma_{c}=\{\alpha:\beta_{0j}(\alpha)=0,~\text{for~any}~j\notin\alpha\}.$$
$\Gamma_{c}$ is the set of all correct models.
However, there could be more than one correct models. Many models  belonging to $\Gamma_{c}$ may also contain some redundant variables $x_{ij}$ that have no effect on response $y_i$. The remaining candidate models can be collected into
$$\Gamma_{w}=\{\alpha:\beta_{0j}(\alpha)\neq0,~\text{for~some}~j\notin\alpha\}.$$
Here $\Gamma_{w}$ is the set of all wrong models, each of which misses at least one non-zero variable $x_{ij}$ that has an effect on $y_i$.

Therefore, it is quite simple to obtain an objective function based on the GLMs likelihood, and use any of the above model selection methods to evaluate the accuracy of all sub-models, unless $p$ can be large enough that implementation of the model selection is computationally prohibitive. In this work, we only study the fixed dimension case, while consistent model selection criteria on the increasing- or high-dimensional case have been studied by \cite{Kim2016}. Now define
\begin{eqnarray}\label{optimal}
{S_n}(\alpha): &=&-{l_n}(\hat \beta (\alpha ))+ C(n,\hat \beta (\alpha ))\nonumber\\
:&=&  \sum\limits_{i = 1}^n {{w_i}} [b(u(x_{i\alpha }^t\hat \beta (\alpha ))) - {y_i}u(x_{i\alpha }^t\hat \beta (\alpha ))] + C(n,\hat \beta (\alpha )),
\end{eqnarray}
where  each $\alpha$ sub-model is estimated via the MLE from the analogue of (\ref{eq:mlog}) for the subvector $\beta(\alpha)$.

Here, the first term in (\ref{optimal}) is the negative log-likelihood that is used to detect the superiority of the model. The second term is a penalty term used to measure the complexity of the model. We assume that $C(n,\hat \beta (\alpha ))$ is an increasing function of the sub-model dimension ${p_\alpha }$. For AIC, $C(n,{{\hat \beta }_n}(\alpha )) = {p_\alpha }$. For BIC, $C(n,{{\hat \beta }_n}(\alpha )) = \frac{1}{2}{p_\alpha }\log n$. For SCC, $C(n,\hat{\beta}_{n}(\alpha))=\log|I_{n}(\hat{\beta}_{n}(\alpha))|/2+\sum\limits_{i=2}^{p_{\alpha}}\log(|\hat{\beta}_{n}(\alpha)_{i}|+\varepsilon n^{-1/4})$ ,where $I_{n}(\beta(\alpha))$ is the Fisher information matrix of $\beta(\alpha)$ and $\varepsilon$ is a specific value to ensure the invariance of the model selection criterion.

Minimize (\ref{optimal}) of all candidate models to get the optimal model, i.e.,
\begin{equation}\label{optimalcri}
\hat \alpha  = \mathop {\arg \min }\limits_{\alpha  \in {\cal A}} {S_n}(\alpha ).
\end{equation}

Under the above criteria, the better the model is fitted, the less complex the model is, i.e., the smaller ${S_n}(\alpha)$ is, the better the model selection effect is. The penalized criterion in (\ref{optimalcri}) shows that the model with good predictability should not only enjoy a small fitting deviation (by the first term in (\ref{optimal})), but also not incorporate unessential variables, i.e., $C(n,\hat{\beta}_{n}(n))$ should be quite small. We know that a model with the smallest ${S_n}(\alpha)$ among all candidate models would be optimal.

Let $\alpha_0$ be the correct model with the smallest size, which is called the \emph{simplest correct model}. That is, the simplest correct model includes exactly all nonzero components of $\beta ({\alpha _0})$. For the simple presentation, we assume that the simplest correct model is unique. Let ${\hat \alpha }$ be the estimate from the model selection criterion (\ref{optimalcri}), then the model selection
procedure is said to be \emph{strongly consistent} if it can select the simplest correct model almost surely, say
\[P\{ \mathop {\lim }\limits_{n \to \infty } \hat \alpha  = {\alpha _0}\}  = 1.\]

The main goal of Section~3 is to evaluate the consistency of the model selection methods that can select the optimal model from all candidate models. Later in the article, some regularity conditions are used to establish model selection criteria and to give some asymptotic results.

\subsection{Regularity Conditions}

We give some general assumptions in advance.
 \begin{itemize}
\item [\textbullet] (H.1): Let ${\lambda _1}\{ S\}  \le  \cdots  \le {\lambda _p}\{ S\} $ be the ordered $p$ eigenvalues of a $p\times p$ symmetric matrix $S$. Then\\
    (i) Let $\{I_{n}(\beta_{0})\}$ be Fisher information of $\beta$ given by \eqref{Fisher-non}. We have
    $\lim\limits_{n\rightarrow\infty}\lambda_{j}\{I_{n}(\beta_{0})\}=\infty$, $j=1,\ldots,p$; (ii) There exist positive constants $d_{1}$,$d_{2}$, that satisfy $d_{1}n\leq\lambda_{p}\{I_{n}(\beta_{0})\}\leq d_{2}n$.

\item [\textbullet] (H.2): Bounded elements of the non-random design: For all $i$, we assume that ${\left\| x_i \right\|_\infty }: = \mathop {\max }\limits_{i,j} \left| {{x_{ij}}} \right| \le L$, where $L$ is a positive constant.

\item [\textbullet] (H.3): We assume that all the coefficients $\beta$ are in parameter space belonging to the space in $\mathbb{R}^p$ such that $\mathop {\sup }\limits_k \left| {\ddot u(x_k^t{\beta})} \right|,\mathop {\sup }\limits_k \left| {\dot u(x_k^t{\beta})} \right|,\mathop {\sup }\limits_k\ddot b(u(x_k^t{\beta })) < \infty $.

\item [\textbullet] (H.4): In the weighted GLMs, if the weights are replaced by some uniformly bounded new weights ${w_k^*}$, $k=1,\cdots,n,$ which are all uniformly bounded, then there exists a positive constant $W$ that satisfies $\mathop {\max }\limits_{1 \le k \le n} w_{k}^{*} \le W$ for all $n$. For notation simplicity, we denote the Fisher information of $\beta$ with old weights in \eqref{Fisher-non} as $I_n^w(\beta ): ={I_n}(\beta )$. We define the weighted Gram matrix as
    \[I_n^{{w^*}}({\beta _0}) = \sum\limits_{k = 1}^n {w_k^*} {{\dot u}^2}(x_k^t{\beta _0})\ddot b(u(x_{k}^{t}\beta_{0})){x_k}x_k^t.\]
The notation $A \prec B$ means that $B-A$ is non-negative definite. The following restricted eigen-value condition is true:
\[{c_l}n{I_p} \prec I_n^{{w^*}}({\beta _0}) \prec {c_u}n{I_p},\]
where $I_p$ is the $p$ dimensional identity matrix, and $c_{l}$, $c_{u}$ are positive constants.

\end{itemize}

Conditions (H.1) below are also proposed in \cite{Wu1999}, \cite{Qian2006} for some special cases of GLMs. Moreover, the combination of (H.1) and (H.2) renders a reasonable convergence rate of estimated coefficients.  The bounded regressors assumption (H.2) is common in some GLM references, see page 46 of \cite{Fahrmeir2001} and Section 2.2.7 of \cite{Chen11} as examples. When we preprocess the raw covariates, ``zero mean'' and ``one variance'' standardizations are required, which evidently and approximately imply the boundedness assumption of covariates. If some predictors have heavy tailed distributions which may be collected in economics and finance, we could take the some-transform such as $f(x)=c\frac{\exp(x)}{1+\exp(x)}$ (where $c$ is a constant) and thus the transformed predictors are approximately seen as bounded variables. For the random design case, we can get all the results by conditioning on $X$.  The (H.3) is important for the consistency results (see \cite{Shao2003}) since the gradient and Hessian for the weighted likelihood in \eqref{eq:mlog} should be bounded to make sure that the optimization problem has a good solution. In the proof, we need to ensure that all the $\mathop {\sup }\limits_k \left| {\ddot u(x_k^t{\beta _0})} \right|, \mathop {\sup }\limits_k \left| {\ddot u(x_k^t{\beta _1})} \right|, \mathop {\sup }\limits_k | {\dot u(x_k^t{\beta _0})} |,\mathop {\sup }\limits_k | {\dot u(x_k^t{\tilde{\beta} _1})} |,\mathop {\sup }\limits_k | {\dot u(x_k^t{\breve{\beta} _1})}|$, and $\mathop {\sup }\limits_k\ddot b(u(x_k^t{\beta _0}))$ are finite, where the values of $\beta_1$ and $\tilde{\beta}_1$ are between the line of $\beta$ and $\beta_0$, and $\beta_1$ is not necessarily equal to $\tilde{\beta}_1$. In addition, the value of $\breve{\beta}_1$ is between the line of $\hat{\beta}$ and $\beta_0$. As for (H.4), with (H.3) it is the theoretical guarantee relating to the proofs of the LIL of the estimator.

\begin{remark}\label{remark 1}
By (H.3) there is a compact subset $K$, not depending on $n$, of the interior of $\Theta$ such that $u(x_k^t{\beta _0})$ belongs to $K$ for all $k$. Let ${e_k} = {y_k} - \dot b(u(x_k^t{\beta _0}))$, the Lemma 6.1 in \cite{Rigollet12} implies
\begin{eqnarray}\label{ei}
\mathop {\sup }\limits_{k \ge 1} E|{e_k}{|^3} = \mathop {\sup }\limits_{k \ge 1} E|{y_k} - \dot b(u(x_k^t{\beta _0})){|^3} \propto \mathop {\sup }\limits_{k \ge 1} {[\ddot b(u(x_k^t{\beta _0}))]^{3/2}} < \infty .
\end{eqnarray}
The moments condition \eqref{ei} is crucial to explain the finite moments behavior of the exponential families. This condition means, in most cases, that there is a compact subset $K$, not depending on $n$, of the interior of $\Theta$ such that $u(x_k^t{\beta _0})$ belongs to $K$ for all $k$. Let $\eta>0$. Condition $\mathop {\sup }\limits_{k \ge 1} E|{e_k}{|^{2 + \eta }} < \infty $ is used in \cite{Yin2006} for establishing the strongly consistent maximum quasi-likelihood estimates.
\end{remark}

\subsection{Independent Observations}
First, we give the following useful strong limit theorems.
\begin{theorem}{(Law of the Iterated Logarithm)}{\label{theorem 1}}
 Assuming conditions (H.1) to (H.4) are satisfied, for any correct model $\alpha\in \Gamma_{c}$,
\begin{eqnarray}\label{3.2}
\|\hat{\beta}(\alpha)-\beta_{0}(\alpha)\|=O(\sqrt{n^{-1}\log\log n}) \qquad  a.s.,
\end{eqnarray}
Furthermore, there is a positive constant $b>0$ such that for arbitrary $\alpha\in \Gamma_{c}$
\begin{eqnarray}\label{3.3}\nonumber
\limsup\limits_{n\rightarrow\infty}\frac{\|\hat{\beta}(\alpha)-\beta_{0}(\alpha)\|}{\sqrt{n^{-1}\log\log n}}=b \qquad a.s.,
\end{eqnarray}
that is, the MLE $\hat{\beta}(\alpha)$ satisfies the LIL.
\end{theorem}

Once Theorem {\ref{theorem 1}} has been constructed, the following almost surely discrepancy between estimated correct sub-model (or incorrect sub-model) log-likelihood and true log-likelihood can be derived as well.

\begin{theorem}\label{theorem 2}
Given the same conditions as Theorem {\ref{theorem 1}}, for any correct model $\alpha\in \Gamma_{c}$,
\begin{eqnarray}\label{3.4}
0 \le {l_n}(\hat \beta (\alpha )) - {l_n}(\beta ({\alpha _0})) = O(\log \log n) \qquad  a.s.,
\end{eqnarray}
where ${l_n}(\hat \beta (\alpha )) = -\sum\limits_{i = 1}^n {{w_i}} [b(u(x_{i\alpha }^t\hat \beta (\alpha ))) - {y_i}u(x_{i\alpha }^t\hat \beta (\alpha ))]$.
\end{theorem}

From the conclusion of Theorem {\ref{theorem 2}}, we can see that the maximum log-likelihood for any correct model $\alpha\in \Gamma_{c}$ is almost surely greater than the unknown log-likelihood of the true model, and the bounds of their differences are almost surely limited by $|O(\log\log n)|$.
\begin{theorem}\label{theorem 3}
Under the same conditions as Theorem {\ref{theorem 1}}, for any incorrect model $\alpha\in \Gamma_{\omega}$, we have
\begin{eqnarray}\label{3.5}
\mathop {\lim \sup }\limits_{n \to \infty } {n^{ - 1}}\{ {l_n}(\hat \beta (\alpha )) - {l_n}({\beta _0})\}  < 0 \qquad  a.s.,
\end{eqnarray}
where ${l_n}({\beta _0}) = -\sum\limits_{i = 1}^n {{w_i}} [b(u(x_i^t{\beta _0})) - {y_i}u(x_i^t{\beta _0})]$.
\end{theorem}

From Theorem {\ref{theorem 3}} we can see that the maximum log-likelihood for all non-correct models $\alpha\in \Gamma_{\omega}$ is almost certainly less than the unknown log-likelihood of the true model. At the same time, when $n$ is sufficiently large, the bounds of their differences will be almost surely limited by $|O(n)|$.
\begin{theorem}\label{theorem 4}
For GLMs with weigted score \eqref{wlog}, under the conditions of Theorem {\ref{theorem 1}}, both BIC and SCC selection criteria are strongly consistent, but AIC selection criterion is not strongly consistent.
\end{theorem}

From Theorem {\ref{theorem 4}} derived by Theorem {\ref{theorem 2}} and Theorem {\ref{theorem 3}}, we know that if we use the penalty log-likelihood method to select the model based on the criterion (\ref{optimalcri}) and the penalty term $C(n,\hat{\beta}(\alpha))$ has an order between $|O(\log\log n)|$ and $|O(n)|$. Moreover, the penalty term $C(n,\hat{\beta}(\alpha))$ is the increasing function of the model dimension, so it will almost surely select the simplest correct model that belongs to $\Gamma_{c}$. The phenomenon that the model selection criteria almost surely choose the simplest correct model is called strong consistency of model selection.

\subsection{Weakly Dependent Observations}
In economics and finance, when responses are collected as time series data with short-term auto-correlation (section 15 in \cite{Hansen18}), heavily relying on the assumption of independent errors are not logical or reasonable. Whereas, most existing references of studying GLMs (includes linear model) are addressed only for independent data, and dependent data are largely not covered. Some work has paid attention to linear models with correlated errors, see \cite{Fan16}. However, in contrast to linear models with dependent responses, the GLMs emphasizing on dependent errors or dependent responses have only been scarcely investigated, see \cite{Kroll2018}. This section focuses on the GLMs endowed with weakly dependent responses, i.e.,
\begin{eqnarray}\label{weakly}
{y_i} = E({y_i}) + {\varepsilon _i},~~E({y_i}) = \dot b(u(x_i^t\beta_0 )),~i = 1,2, \cdots ,n,
\end{eqnarray}
where ${\varepsilon _i}$ are weakly dependent error sequence with zero mean and $\{ {x_i}\} $ are fixed.

Although, the weakly dependence cannot be solved by MLE approach due to the dependent structure, the likelihood does not enjoy the form of production. However, by applying quasi-likelihood method, we can still be able to get a desired estimator which is consistent under some regularity conditions. The quasi-likelihood estimator $\hat{\beta}_{\rm Q}(\alpha)$ is the solution to follow estimating equation:
 \begin{equation}\label{eq:likelihood}
S_{n}(\beta_0):=\sum\limits_{i=1}^{n}w_{i}x_{i}\dot{u}(x_{i}^{t}\beta_0)[y_{i}-E(y_{i})]=0.
 \end{equation}
The inference from the solution of \eqref{eq:likelihood} is a typical example and is the maximum likelihood estimate for independent responses. More details about quasi-likelihood can be seen in \cite{Fahrmeir2001}, \cite{Chen11}.

Next, we still want to apply LIL that would imply model selection consistency by using the similar approach mentioned in the case for independent observations. In order to present the dependence structure, we point out some notations and definitions for measuring dependence.

Let $(\Omega,\mathcal{F},P)$ be a probability space, and $\mathcal{B}$, $\mathcal{C}$ two sub-$\sigma$-fields of $\mathcal{F}$. Let ${L_2}({\cal B})$ be a set of all ${\cal B}$-measurable random variables with finite 2nd moments.

We introduce the notation of $\rho$-mixing and $\alpha$-mixing coefficients (strong mixing coefficient):
\begin{eqnarray}\nonumber
\rho ({\cal B},{\cal C}) = \mathop {\sup }\limits_{X \in {L_2}({\cal B}),Y \in {L_2}({\cal C})} \frac{{\left| {E(XY) - E(X)E(Y)} \right|}}{{\sqrt {Var(X)Var(Y)} }},\\
\alpha ({\cal B},{\cal C}) = \mathop {\sup }\limits_{B \in {\cal B},C \in {\cal C}} \left| {P(B \cap C) - P(B)P(C)} \right|.
\end{eqnarray}

Let $\mathcal{F}_a^b = \sigma ({X_t},a \le t \le b)$ and $\mathbb{N}$ be the set of all positive integers. The random sequence $\{ {X_t},t \in \mathbb{N}\} $ is said to be $\alpha$-\emph{mixing} or \emph{strongly mixing}, if
\begin{eqnarray}\nonumber
\alpha(n): = \mathop {\sup }\limits_{k \in \mathbb{N}} \alpha (\mathcal{F}_1^k,\mathcal{F}_{k + n}^\infty ) \to 0,~~n \to \infty .
\end{eqnarray}
Next, a sequence of random variables $\{ {X_t},t \in \mathbb{N}\} $ is called a $\rho$-\emph{mixing process} if
\begin{eqnarray}\nonumber
\rho (n): = \mathop {\sup }\limits_{k \in \mathbb{N}} \rho (\mathcal{F}_1^k,\mathcal{F}_{k + n}^\infty ) \to 0,~~n \to \infty .
\end{eqnarray}

A sequence of strictly stationary random variables $\{ {X_t},t \in \mathbb{N}\} $ are said to be $m$-dependent, if for any two subsets $B, C \subset \mathbb{N}$ and $\mathop {\inf }\limits_{{t_1} \in B,{t_2} \in C} \left| {{t_1} - {t_2}} \right| > m$, the process $\{ {X_{t_1}},{t_1} \in B\} $ and $\{ {X_{t_2}},{t_2} \in C\} $ are independent.

For detailed theories and examples (such linear time series) about $\alpha$-mixing process, $\rho$-mixing process and other mixing processes, we refer reading to \cite{Lin1997} and \cite{Bosq1998}. Notice that, we have $\alpha ({\cal A},{\cal B}) \le \frac{{\rm{1}}}{{\rm{4}}}\rho ({\cal A},{\cal B})$ and then the $\rho$-mixing process implies the $\alpha$-mixing process. For simplicity, we would restrict our study to strictly stationary sequences in this section.

To prepare for the asymptotic theories for the weakly dependent case, we need some additional regularity assumptions for our technique proofs.
\begin{itemize}

 \item [\textbullet] (H.5): The ${y _i}$ satisfies the $\rho$-mixing condition with geometric decay: $\rho(n) =O({r^{ - n}})$.

  \item [\textbullet] (H.6): The ${y_i}$ is $m$-dependent.
\end{itemize}

Since weak dependency makes the problem more complex, we assume that the design matrix $X$ is viewed as
a non-random matrix. A similar assumption is given by \cite{Fan16}. Differ from Condition 3 in \cite{Fan16}, they assumed that the $\alpha$-mixing condition is by exponential decay $\alpha (m) = O({e^{ - a{m^b}}})$ where $a, b$ are positive constant. Our geometric rate of decay is not sharper than exponential decay, which means that the dependence of our response is allowed to be stronger. Property (H.5) or (H.6) holds for the $\varepsilon_{i}$ if and only if it holds for the $y_{i}$. The same assertion holds for (H.6). It should be note that for (H.6) if an $m$-dependent sequence has $\rho (n) = O({r^{ - n}})$ for $n > m$, then the sequence is $\rho$  mixing. But it may not be $\rho$-mixing with geometric decay $\rho(n) =O({r^{ - n}})$. The (H.6) implies the $\rho$-mixing with truncated decay: $\rho (n) = O({a_n} \cdot 1\{ n \le m\} )$ for some sequence $\{ {a_n}\}$.

Now, we present the main results below which are the same as situation of independent observations.

\begin{theorem}{(Law of the Iterated Logarithm)}\label{theorem 5}
 Assuming conditions (H.1) to (H.4) are satisfied. For weakly dependent GLMs with estimating equation \eqref{eq:likelihood}, the $\rho$-mixing responses $\{ {y_i}\} _{i = 1}^n$ with additional requirement (H.5) or the $m$-dependent responses $\{ {y_i}\} _{i = 1}^n$ with additional requirement (H.6),  we have
\begin{eqnarray}\label{3.2wd}
\|\hat{\beta}(\alpha)-\beta_{0}(\alpha)\|=O(\sqrt{n^{-1}\log\log n}) \qquad  a.s.
\end{eqnarray}
for any correct model $\alpha\in \Gamma_{c}$.

Furthermore, there is a positive constant $b>0$ such that for any arbitrary $\alpha\in \Gamma_{c}$
\begin{eqnarray}\label{3.3wd}\nonumber
\limsup\limits_{n\rightarrow\infty}\frac{\|\hat{\beta}(\alpha)-\beta_{0}(\alpha)\|}{\sqrt{n^{-1}\log\log n}}=b \qquad a.s.,
\end{eqnarray}
where MLE estimates $\hat{\beta}(\alpha)$ satisfy the law of iterated logarithm.
\end{theorem}
Once the LIL for weakly dependent GLMs is verified, the following strong convergence for the BIC or SCC model selection is similarly to be obtained, by checking the Theorem {\ref{theorem 2}} and Theorem {\ref{theorem 3}}.

\begin{theorem}\label{theorem 6}
In the weakly dependent GLMs with estimating equation \eqref{eq:likelihood}, considering the $\alpha$-mixing responses or $m$-dependent $\{ {y_i}\} _{i = 1}^n$ error sequences. Under the same conditions of Theorem {\ref{theorem 5}}, both BIC and SCC criteria are strongly convergent for the model selection, but not strongly convergent for the AIC criterion.
\end{theorem}

\section{SIMULATION STUDY}
 The purpose for this section is to examine the difference between the performance of BIC and AIC in the variable selection. The BIC criterion for MLE $\hat \beta^{nb} (\alpha)$ of NBR is defined as
\begin{align*}
{\rm{BIC}}[\hat \beta^{nb} (\alpha)] :&= - \frac{1}{n}\sum\limits_{i = 1}^n {[{y_i}x_i^t\hat \beta ^{nb} (\alpha)  - (\theta  + {Y_i})\log (\theta  + {e^{x_i^t\hat \beta^{nb} (\alpha) }})]}+ \frac{{\log n}}{n}\hat{\rm{df}}[\hat \beta^{nb} (\alpha)],
\end{align*}
where $\hat{\rm{df}}(\hat \beta^{nb} (\alpha)):=||\hat \beta^{nb} (\alpha)|{|_0}$ is the number of coefficients in the model. Similarly, the BIC criterion for MLE $\hat \beta^{pb} (\alpha)$ of probit model is given by
\begin{align*}
{\rm{BIC}}[\hat \beta^{pb} (\alpha)]:&= - \frac{1}{n}\sum\limits_{i = 1}^n {\{ {y_i}\log \Phi (x_i^t\hat \beta^{pb} (\alpha) ) + (1 - {y_i})\log [1 - \Phi (x_i^t\hat \beta^{pb} (\alpha))]\} }\\
&+ \frac{{\log n}}{n}\hat{\rm{df}}[\hat \beta^{pb} (\alpha)]
\end{align*}
where $\Phi(t)$ is the cumulative distribution function of $N(0,1)$.

For weakly dependent responses, we consider dependent linear model with $MR(p)$ time series ($p=2,3$, see section 15.4 \cite{Hansen18}) as the error sequence $ \{ {\varepsilon _i}\}$, that is
\begin{eqnarray}\label{weaklyLM}
{y_i} = x_i^t\beta  + {\varepsilon _i},~i = 1,2, \cdots ,n,
\end{eqnarray}
We solve \eqref{weaklyLM} (denote $\hat \beta^{dl} (\alpha)$) by the OLS method which is equivalent to the estimation from quasi-likelihood. The BIC criterion for dependent linear model is
\[{\rm{BIC}}[\hat \beta^{dl} (\alpha)]:= \log [\frac{1}{n}\sum\limits_{i = 1}^n {{{({y_i} - x_i^t\hat \beta (\alpha ))}^2}} ] + \frac{{\log n}}{n}\hat{\rm{df}}[\hat \beta^{dl} (\alpha)]\]

At the same time, for comparison,  we consider the AIC criterion which replaces the BIC penalty function $\frac{{\log n}}{n}\hat{\rm{df}}[\hat \beta (\alpha)]$ by the penalty function $\frac{1}{n}\hat{\rm{df}}[\hat \beta (\alpha)]$.

In the simulation, we add three more redundant variables $X_{4i}, X_{5i}, X_{6i}$ in each observation, which are also independent and distributed $U(0,1)$. The negative binomial regression model and the probit regression model are considered in our simulation and we set the $\mathbf{\beta}$ in the negative binomial regression and the probit model to be 0.5. The $\theta$ in the NBR is 10 in the simulation. As for dependent linear model, we consider the error term which satisfies $MR(2)$ and $MR(3)$ with parameters (0.5, 0.3) and (0.5, 0.3, 0.2).

We run the simulation for 500 times. Each time all ($2^6-1=63$) kinds of combinations of the variables denote a model selection to be ``correct'' if
the method successfully selects exactly $\{X_{i1}, X_{i2}, X_{i3}\}$. If the criterion picks other variables besides $\{X_{i1}, X_{i2}, X_{i3}\}$, we say it is an ``overfit''. The rest is denoted as ``underfit''. The MSE is also calculated in our simulation. From table \ref{tab:variableselection} we can see that, in comparison with AIC, BIC criterion always performs better in selecting the true model. At the same time, with the growth of sample size $N$, the  ``correct'' of BIC criterion gradually converges to 1.

\begin{table}[!htbp]
\centering
\caption{Simualtion results for variable selection.}\label{tab:variableselection}
\begin{tabular}{ccccccc}
\toprule
Model & Method&Sample Size& Correct&Overfit&Underfit&MSE\\
\midrule
\hline
\multirow{4}{*}{NBR ($\beta=0.5,\theta=10$)}& \multirow{2}{*}{BIC}& 100& 0.8920&0.1053&0.0037&0.0164\\
 &  & 300& 0.9427&0.0573&0.0000&0.0064\\
& \multirow{2}{*}{AIC}& 100& 0.5617&0.4373&0.0010&0.0287\\
 &  & 300& 0.5863&0.4137&0.0000&0.0197\\
\hline
\multirow{4}{*}{Probit model($\beta=0.5$)}& \multirow{2}{*}{BIC}&100 &0.7106&0.0816&0.2078&0.0190\\
 &  & 300&0.9444&0.0000&0.0556&0.0094\\
& \multirow{2}{*}{AIC}&100 &0.5420&0.4164&0.0416&0.0155\\
 &  & 300&0.5928&0.0000&0.4072&0.0193\\
 \hline
\multirow{4}{*}{Dependent LM (MR(2) error)}& \multirow{2}{*}{BIC}&100 &0.8460&0.1002&0.0538&0.0178\\
 &  & 300&0.9468&0.0532&0.0000&0.0104\\
& \multirow{2}{*}{AIC}&100 &0.5590&0.4330&0.0080&0.0222\\
 &  & 300&0.5928&0.0000&0.4072&0.0193\\
 \hline
\multirow{4}{*}{Dependent LM (MR(3) error)}& \multirow{2}{*}{BIC}&100 &0.8366&0.1024&0.0610&0.0181\\
 &  & 300&0.9510&0.0490&0.0000&0.0110\\
& \multirow{2}{*}{AIC}&100 &0.5642&0.4254&0.0181&0.0254\\
 &  & 300&0.5890&0.4110&0.0000&0.0116\\
 \hline
\end{tabular}
\end{table}

Since the AIC/BIC analysis of GLMs have been applied widely to real data, we ignore the real data analysis in this work, see \cite{Tutz2011}.

\section{DISCUSSION}

This paper presents an unified theory for studying strong limit theorems for independent and dependent GLMs responses and we obtain the LIL for the maximum likelihood estimate of the regression coefficients.
The convergence rate of regression coefficients estimators is shown to be $O(\sqrt{n^{-1}\log\log n})$ almost surely. This LIL is employed to establish the strong consistency of BIC and SCC selection criteria.

Using the same techniques in this paper, we can also prove that the BIC procedure for sub-sampling based GLMs estimates in \cite{Ai2018} also enjoys the strong consistency of model selection. Further research of this work may be the LIL for the maximum likelihood estimator in errors-in-variables GLMs under mild conditions, where \cite{Miao2011} only considered the case of simple linear regression.

\section{Acknowledgements}
The authors would to show their hearty thanks to Professor YueHua Wu for her helpful comments and mentioned references \cite{Wu1999} and \cite{Rao1989} which substantially improve this work.

\section{APPENDIX}

\begin{appendix}
The key of our proofs is the mathematical analysis of the negative likelihood function, and the local quadratic approximation is frequently used. The main difficulty lies in how to standardize the score function and how to establish a consistent bound on the the difference of two log-likelihood functions almost exactly. The techniques of using some concentration inequalities and the error bound of log-likelihood ratio are broadly employed in literatures, see \cite{Rao1989}, \cite{Rao1992}, \cite{Wu1999}, compare \cite{Qian2006} for similar considerations.

First, we define a positive sequence $\{\tau_{n}\}$ such that
\begin{eqnarray}\label{twosequence}
\tau_{n}\uparrow \infty~\text{(or~even~a large constant)}, \quad \tau_{n}(n^{-1}\log\log n)^{\frac{1}{2}}\downarrow 0.
\end{eqnarray}
In the sequel, we introduce a sequence of $L_2$ ball with radius being proportional to $\{\tau_{n}\}$:
\[{B_n} = \{ \beta :\|\beta  - {\beta _0}\| \le {\tau _n}{({n^{ - 1}}\log \log n)^{\frac{1}{2}}}\} \]
and denote $\partial B_{n}=\{\beta:\|\beta-\beta_{0}\|=\tau_{n}(n^{-1}\log\log n)^{\frac{1}{2}}\}$ as its boundary. Obviously, we have $B_{1}\supset B_{2}\supset B_{3}\supset\cdots\supset B_{n}$.

Further, we define the log-likelihood ratio of two log-likelihoods in terms of $\beta$ and $\beta_{0}$:
\begin{eqnarray}\label{3.7}\nonumber
{K_n}(\beta ,{\beta _0}):=l_{n}(\beta_{0})-l_{n}(\beta)=\sum\limits_{k=1}^{n}w_{k}\{b(u(x_{k}^{t}\beta))-b(u(x_{k}^{t}\beta_{0}))-y_{k}[u(x_{k}^{t}\beta)-u(x_{k}^{t}\beta_{0})]\}.\nonumber
\end{eqnarray}
Let $K(t,s)=b(t)-b(s)-\dot{b}(s)(t-s)$, then
\begin{eqnarray}\nonumber
K(u(x_k^t\beta ),u(x_k^t{\beta _0})) = b(u(x_k^t\beta )) - b(u(x_k^t{\beta _0})) - \dot b(u(x_k^t{\beta _0}))[u(x_k^t\beta ) - u(x_k^t{\beta _0})].
\end{eqnarray}
So ${K_n}(\beta ,{\beta _0})$ can be rewritten as
\begin{eqnarray}\nonumber
{K_n}(\beta ,{\beta _0})&=&\sum\limits_{k=1}^{n}w_{k}\{K(u(x_{k}^{t}\beta),u(x_{k}^{t}\beta_{0}))
+(\dot{b}(u(x_{k}^{t}\beta_{0}))-y_{k})[u(x_{k}^{t}\beta)-u(x_{k}^{t}\beta_{0})]\}\nonumber\\
&=&\sum\limits_{k=1}^{n}w_{k}K(u(x_{k}^{t}\beta),u(x_{k}^{t}\beta_{0}))
-\sum\limits_{k=1}^{n}w_{k}(y_{k}-\dot{b}(u(x_{k}^{t}\beta_{0})))[u(x_{k}^{t}\beta)-u(x_{k}^{t}\beta_{0})]\nonumber
\end{eqnarray}
Our main idea is to assume that the $u(\cdot)$ function is sufficiently smooth, then we can use Taylor expansion for $u(t)$ at $t=s$:
\begin{eqnarray}\nonumber
u(t)=u(s)+\dot{u}(s)(t-s)+\frac{1}{2}\ddot{u}(\zeta)(t-s)^2,
\end{eqnarray}
where $\zeta$ is between $t$ and $s$.

Let $t=x_{k}^{t}\beta$, $s=x_{k}^{t}\beta_{0}$, $\zeta=x_{k}^{t}\beta_{1}$.  The value of $\beta_{1}$ is between $\beta$ and $\beta_{0}$. Then by Taylor's expansion we have
\begin{eqnarray}\label{3.7}\nonumber
{K_n}(\beta ,{\beta _0})&=&\sum\limits_{k = 1}^n {{w_k}} K(u(x_k^t\beta ),u(x_k^t{\beta _0})) - \sum\limits_{k = 1}^n {{w_k}} ({y_k} - \dot b(u(x_k^t{\beta _0})))[\dot u(x_k^t{\beta _0})(x_k^t\beta  - x_k^t{\beta _0})]\nonumber \\
&-& \frac{1}{2}\sum\limits_{k = 1}^n {{w_k}} \ddot u(x_k^t{\beta _1})[{y_k} - \dot b(u(x_k^t{\beta _0}))]{(x_k^t\beta  - x_k^t{\beta _0})^2}\nonumber\\
&=:&{R_{1n}}(\beta ) + {R_{2n}}(\beta ) + {R_{3n}}(\beta ),
\end{eqnarray}
where
\[\begin{array}{l}
{R_{1n}}(\beta ) = \sum\limits_{k = 1}^n {{w_k}} K(u(x_k^t\beta ),u(x_k^t{\beta _0})),\\
{R_{2n}}(\beta ) =  - \sum\limits_{k = 1}^n {{w_k}} \dot u(x_k^t{\beta _0})[{y_k} - \dot b(u(x_k^t{\beta _0}))]x_k^t(\beta  - {\beta _0}),\\
{R_{3n}}(\beta ) = -\frac{1}{2}\sum\limits_{k = 1}^n {{w_k}} \ddot u(x_k^t{\beta _1})[{y_k} - \dot b(u(x_k^t{\beta _0}))]{(x_k^t\beta  - x_k^t{\beta _0})^2}.
\end{array}\]

\section{Lemmas}

Before presenting the proof of the main results, we give some lemmas which will be useful in the following sections.

The following lemmas that we will apply is a result in convex optimizations.

\begin{lemma}\label{lemma a.1}{(Lemma 3.3.10 of \cite{Wright2017})}
Assume that function $f(\cdot)$ is strongly convex and $\dot{f}(\cdot)$ is Lipschitz
continuously differentiable, then $f(x,y)$ has the following upper and lower bounds for arbitrary vectors x, y,
\begin{eqnarray}\label{3.8}
\frac{1}{2}C_{l}\|y-x\|^{2} \leqslant f(y)-f(x)-\dot{f}(x)(y-x) \leqslant \frac{1}{2}C_{u}\|y-x\|^{2},
\end{eqnarray}
where $C_{u}$, $C_{l}$ are positive constants.
\end{lemma}

The following Petrov's law of iterated logarithm for sums of independent random variables can be found in  p274 of \cite{Stout1974}.
\begin{lemma}{(Petrov's LIL)}\label{lemma a.2}
Let $\left( X _ { i } \right) _ { i \geqslant 1 }$ be a sequence of centered independent random variables and define ${S_n} = \sum\nolimits_{i = 1}^n {{X_i}}$. Suppose that:\\

(i) $\mathop {\lim }\limits_{n \to \infty } \frac{1}{n}{\rm{Var}}{S_n} = {\sigma ^2} < \infty $; (ii) $\mathop {\sup }\limits_{i \ge 1} E|X_i|^{2 + \eta } < \infty $ for some positive $\eta $.\\
 Then
 \[\mathop {\limsup }\limits_{n \to \infty } \frac{{{S_n}}}{{\sqrt {2n\log \log n} }} = \sigma~~\text{a.s.}\]
\end{lemma}


\begin{corollary}\label{cor 1}
Under (H.1)-(H.4), we have for $j=1,2,\cdots,p.$
\begin{eqnarray}\label{3.11}
\mathop {\lim \sup }\limits_{n \to \infty } \frac{{| {\sum\limits_{k = 1}^n {{w_k}} \dot u(x_k^t{\beta _0})[{y_k} - \dot b(u(x_k^t{\beta _0}))]{x_{kj}}} |}}{{\sqrt {2{I_n}({\beta _0})(j,j)\log \log {I_n}({\beta _0})(j,j)} }} = 1 \qquad a.s.
\end{eqnarray}
where $w_{k}$ is the $k$th component of weight $w$, $x_{kj}$ is the $j$th element of vector $x_{k}$, and $I_{n}(\beta_{0})(j,j)$ is the $(j,j)$th element of the Fisher information matrix $I_{n}(\beta_{0})$. The above formula shows that $\{w_{k}\dot{u}(x_{k}^{t}\beta_{0})(y_{k}-\dot{b}(u(x_{k}^{t}\beta_{0})))\}$ meets the LIL, so we have
\begin{eqnarray}\label{3.12}
S_{n}(\beta_{0})&=&\sum\limits_{k=1}^{n}w_{k}\dot{u}(x_{k}^{t}\beta_{0})(y_{k}-\dot{b}(u(x_{k}^{t}\beta_{0})))x_{k}=O(\sqrt{n\log\log n}){1_p}~a.s.
\end{eqnarray}
where ${1_p} := {(1, \cdots ,1)^t}$.
\end{corollary}

\begin{corollary}\label{cor 2}
Under (H.1)-(H.4), we have
\begin{eqnarray}\label{3.12}
\mathop {\lim \sup }\limits_{n \to \infty } \frac{{|\sum\limits_{k = 1}^n {{w_k}} \ddot u(x_k^t{\beta _1})[{y_k} - \dot b(u(x_k^t{\beta _0}))]||x_k^t||_1^2|}}{{\sqrt {2n\log \log n} }} <\infty.
\end{eqnarray}
\end{corollary}

\section{Proof of Corollary {\ref{cor 1}}}

In the following, we only need to verify (\ref{3.11}). Note that the mean and variance of $\{y_{k}\}_{k=1}^{n}$ are $E(y_{k})=\dot{b}(u(x_{k}^{t}\beta_{0})), Var(y_{k})=\ddot{b}(u(x_{k}^{t}\beta_{0}))$ and the definition of $I_{n}(\beta_{0})$, for $j=1,\cdots,p$, let
\begin{eqnarray}\label{3.13}
{S_{nj}}: = \sum\limits_{k = 1}^n {{w_k}} \dot u(x_k^t{\beta _0})[{y_k} - \dot b(u(x_k^t{\beta _0}))]{x_{kj}}
\end{eqnarray}

According to condition (H.4) and (H.1)(i), we have
\begin{eqnarray}\label{3.14}\nonumber
\frac{1}{n}{\rm{Var}}({S_{nj}}) &=&\frac{1}{n}\sum\limits_{k=1}^{n}w_{k}^{2}\dot{u}(x_{k}^{t}\beta_{0})^{2}E(y_{k}-E(y_{k}))^{2}x_{kj}^{2}\nonumber\\
&\leq & \frac{1}{n}W\sum\limits_{k=1}^{n}w_{k}\dot{u}(x_{k}^{t}\beta_{0})^{2}\ddot{b}(u(x_{k}^{t}\beta_{0}))x_{kj}^{2}=\frac{1}{n}WI_{n}(\beta_{0})(j,j)=O(1).
\end{eqnarray}
Thus condition (i) of Lemma {\ref{lemma a.2}} is verified.

To check (ii) of Lemma {\ref{lemma a.2}}, we only need to use Remark {\ref{remark 1}} in the Section 3.2.


\section{Proof of Corollary {\ref{cor 2}}}

The result (\ref{3.12}) comes from (\ref{3.11}) and condition (H.1). The proof is similar to Corollary {\ref{cor 1}}.

\section{Proof of Theorem {\ref{theorem 1}}}
Judging from the lemmas above, now we are able to prove Theorem {\ref{theorem 1}}. For notation simplicity, we write $\beta (\alpha ),{x_{k\alpha }}$ as $\beta ,{x_k}$, respectively.

Firstly, we give the proof of the non-canonical link GLMs.

Let $t=u(x_{k}^{t}\beta)$, $s=u(x_{k}^{t}\beta_{0})$, ${\tilde \beta _1}$ be the value between $\beta_{0}$ and $\beta$. Then it can be seen from Lemma {\ref{lemma a.1}} that
\begin{eqnarray}\label{3.24}
K(u(x_{k}^{t}\beta),u(x_{k}^{t}\beta_{0}))&:=&b(u(x_{k}^{t}\beta))-b(u(x_{k}^{t}\beta_{0}))-\dot{b}(u(x_{k}^{t}\beta_{0}))(u(x_{k}^{t}\beta)-u(x_{k}^{t}\beta_{0}))\nonumber\\
&\geq& \frac{1}{2}C_l(u(x_{k}^{t}\beta)-u(x_{k}^{t}\beta_{0}))^{2}=frac{1}{2}{C_l}{{\dot u}^2}(x_k^t{\tilde \beta_1})(x_k^t\beta  - x_k^t{\beta _0}{)^2},
\end{eqnarray}
where the last equality is from the Taylor's expansion.

Then
\begin{eqnarray}\label{3.26}\nonumber
R_{1n}(\beta)I(\beta\in\partial B_{n})&=&\sum\limits_{k=1}^{n}w_{k}K(u(x_{k}^{t}\beta),u(x_{k}^{t}\beta_{0}))I(\beta\in\partial B_{n})\nonumber\\
&\geq&\frac{1}{2}C_{l}\sum\limits_{k=1}^{n}w_{k}{{\dot u}^2}(x_k^t{\tilde{\beta} _1})(x_{k}^{t}\beta-x_{k}^{t}\beta_{0})^{2}I(\beta\in\partial B_{n})\nonumber\\
&=&O(1){(\beta  - {\beta _0})^t}[\sum\limits_{k = 1}^n {w_k^{}} {{\dot u}^2}(x_k^t{\tilde{\beta} _1})x_k^{}x_k^t](\beta  - {\beta _0})I(\beta  \in \partial {B_n}).\nonumber\\
~\label{3.26}
\end{eqnarray}
Denote a weighted Gram matrix as
\begin{eqnarray*}
I_n^{{w^*}}({\beta _0}) := \sum\limits_{k = 1}^n {w_k^*} {{\dot u}^2}(x_k^t{\beta _0})\ddot b(u(x_{k}^{t}\beta_{0})){x_k}x_k^t.
\end{eqnarray*}

Using (H.3) and (H.4), we have ${c_l}nI_{p} \prec I_n^{{w^*}}({\beta _0})$ for large $n$. Therefore,
\begin{eqnarray}\label{3.26b}\nonumber
R_{1n}(\beta)I(\beta\in\partial B_{n})&=&O(1){(\beta  - {\beta _0})^t}I_n^{{w^*}}({\beta _0})(\beta  - {\beta _0})I(\beta  \in \partial {B_n})\nonumber\\
&\geq& O(n)\|\beta-\beta_{0}\|^{2}I(\beta\in\partial B_{n}).
\end{eqnarray}
By $\partial B_{n}=\{\beta:\|\beta-\beta_{0}\|=\tau_{n}\sqrt{n^{-1}\log\log n}\}$, then
\begin{eqnarray}\label{3.27}
{R_{1n}}(\beta )I(\beta  \in \partial {B_n}) \ge O(1)\tau _n^2\log \log n.
\end{eqnarray}
For ${R_{2n}}(\beta )$, according to the Cauchy's inequality and the LIL, we have
\begin{eqnarray}\label{3.28}\nonumber
|{R_{2n}}(\beta )|I(\beta  \in \partial {B_n}) &\le& \| {\sum\limits_{k = 1}^n {{w_k}} \dot u(x_k^t{\beta _0})({y_k} - \dot b(u(x_k^t{\beta _0}))){x_k}} \| \cdot \| {\beta  - {\beta _0}} \|I(\beta  \in \partial {B_n})\nonumber\\
&=&O(\sqrt{n\log\log n})\tau_{n}\sqrt{n^{-1}\log\log n}=O(\tau_{n})\log\log n~a.s..
\end{eqnarray}
Then as for ${R_{3n}}(\beta )$, by the Cauchy's inequality again, it gives
\begin{eqnarray}\nonumber
|{R_{3n}}(\beta )|I(\beta  \in \partial {B_n}) &=& |\frac{1}{2}\sum\limits_{k = 1}^n {{w_k}} \ddot u(x_k^t{\beta _1})[{y_k} - \dot b(u(x_k^t{\beta _0}))]{(x_k^t\beta  - x_k^t{\beta _0})^2}|I(\beta  \in \partial {B_n})\nonumber\\
&\le& |\frac{1}{2}\sum\limits_{k = 1}^n {{w_k}} \ddot u(x_k^t{\beta _1})[{y_k} - \dot b(u(x_k^t{\beta _0}))]| \cdot |x_k^t(\beta  - {\beta _0}){|^2} I(\beta  \in \partial {B_n}) \nonumber\\
&=&|\frac{1}{2}\sum\limits_{k = 1}^n {{w_k}} \ddot u(x_k^t{\beta _1})[{y_k} - \dot b(u(x_k^t{\beta _0}))]||x_k^t||_1^2| \cdot ||\beta  - {\beta _0}||_\infty ^2. \label{R3S}
\end{eqnarray}
By Corollary {\ref{cor 2}}, we get
\begin{eqnarray}\label{LIL1}
| \sum\limits_{k = 1}^n {{w_k}} \ddot u(x_k^t{\beta _1})[{y_k} - \dot b(u(x_k^t{\beta _0}))]||x_k^t||_1^2 | = O(\sqrt {n\log \log n} ).
\end{eqnarray}
So the last equality in \eqref{R3S} is obtained.

Therefore,
\begin{eqnarray}\label{3.281}
|{R_{3n}}(\beta )|I(\beta  \in \partial {B_n}) &=& O(\sqrt {n\log \log n} ) \cdot {\left\{ {{\tau _n}\sqrt {{n^{ - 1}}\log \log n} } \right\}^2}\nonumber\\
&=& O(1) n^{-0.5}\cdot \tau _n^2{(\log \log n)^{1.5}} \qquad a.s..
\end{eqnarray}
So, for large $n$, we get
\begin{align}\label{3.32}
{K_n}(\beta ,{\beta _0})I(\beta  \in \partial {B_n}) &= \{ {R_{1n}}(\beta ) + {R_{2n}}(\beta ) + {R_{3n}}(\beta )\} I(\beta  \in \partial {B_n})\nonumber\\
 &\ge O(\tau _n^2){\log \log n}~~a.s..
\end{align}
since $\tau _n^2\log \log n >  > \tau _n^{}\log \log n$ and $\tau _n^{} \to \infty ~(\text{or~some~big~constant})$.

From Lemma {\ref{lemma a.1}}, we see that ${K_n}(\beta ,{\beta _0})$ is convex and then ${K_n}(\beta_0 ,{\beta _0})=0$. The (\ref{3.32}) shows that the ${\hat \beta }$ minimizes $\frac{1}{n}{K_n}(\beta ,{\beta _0})$ to be contained in a subset of $B_{n}$ almost surely. Then
\begin{eqnarray}\label{3.33}
\|\hat{\beta}-\beta_{0}\|\leq\tau_{n}\sqrt{n^{-1}\log\log n} \qquad a.s.
\end{eqnarray}
Since the sequence $\{\tau_{n}\}$ is chosen to be as slow as possible and tends to infinity (or tends to a constant $C$ such that $O(1)\tau _n^2\log \log n + O(1)\tau _n^{}\log \log n + O(1)n^{-0.5}\tau _n^{2}{(\log \log n)^{1.5}} \ge {C^2}{\log \log n}$ and $C$ is depended by the constants on \eqref{3.27}, \eqref{3.28} and \eqref{3.281}, then (\ref{3.33}) can be used to get the establishment of (\ref{3.2}).

Furthermore, we need to prove
\begin{eqnarray}\label{3.34}
\limsup\limits_{n\rightarrow\infty}\frac{\|\hat{\beta}(\alpha)-\beta_{0}(\alpha)\|}{\sqrt{n^{-1}\log\log n}}=d \ne 0  \qquad a.s.
\end{eqnarray}
The following is obtained by contradiction. If it is assumed that the above formula is not established, then
\begin{eqnarray}\label{3.35}
\|\hat{\beta}-\beta_{0}\|=o(\sqrt{n^{-1}\log\log n}) \qquad a.s.
\end{eqnarray}
By applying the result of Lemma {\ref{lemma a.1}}, let $\breve{\beta} _{1}$ be the value between $\beta_{0}$ and $\hat{\beta}$, we have
\begin{eqnarray}\label{3.36}\nonumber
{R_{1n}}(\hat{\beta})&\leq&\frac{1}{2}C_{u}\sum\limits_{k=1}^{n}w_{k}{{\dot u}^2}(x_k^t{\breve{\beta} _1})(x_{k}^{t}\hat{\beta}-x_{k}^{t}\beta_{0})^{2}\nonumber\\
&=&\frac{1}{2}C_{u}(\hat{\beta}-\beta_{0})^{t}\sum\limits_{k=1}^{n}w_{k}\dot{u}^{2}(x_{k}^{t}\breve{\beta}_{1})x_{k}x_{k}^{t}(\hat{\beta}-\beta_{0})\nonumber\\
&=&O(n)\|\hat{\beta}-\beta_{0}\|^{2}= o(1)\cdot\log\log n \qquad a.s.
\end{eqnarray}
where the last equality stems from \eqref{3.35}.

Similarly, the order of ${R_{2n}}(\hat{\beta}), {R_{3n}}(\hat{\beta})$ can also be drawn
\begin{eqnarray}\label{3.38}
{R_{2n}}(\hat{\beta})&=&O(\sqrt{n\log\log n})\cdot o(\sqrt{n^{-1}\log\log n})=o(1)\cdot\log\log n ~a.s.\\
{R_{3n}}(\hat{\beta})&=&o(1)\cdot\log\log n ~ a.s.
\end{eqnarray}

So under the assumption of (\ref{3.35}),
\begin{eqnarray}\label{smallo}
K_{n}(\hat{\beta},\beta_{0})=R_{1n}(\hat{\beta})+R_{2n}(\hat{\beta})+R_{3n}(\hat{\beta})=o(1)\log\log n \qquad a.s.
\end{eqnarray}
On the other hand, from the formula (\ref{3.11}) in Corollary {\ref{cor 1}}, we know that there exists a sub-sequence $\{n_{i}\uparrow\infty\}$ satisfying
\begin{eqnarray}\label{LIL2}
\lim\limits_{i\rightarrow\infty}\frac{|\sum\limits_{k=1}^{n_{i}}w_{k}\dot{u}(x_{k}^{t}\beta_{0})(y_{k}-\dot{b}(u(x_{k}^{t}\beta_{0})))x_{k1}|}{\sqrt{2I_{n_{i}}(\beta_{0})(1,1)\log\log I_{n_{i}}(\beta_{0})(1,1)}}=1  \qquad a.s.
\end{eqnarray}
So when $n_{i}$ is big enough, we have
\begin{eqnarray}\label{3.42}
\frac{|\sum\limits_{k=1}^{n_{i}}w_{k}\dot{u}(x_{k}^{t}\beta_{0})(y_{k}-\dot{b}(u(x_{k}^{t}\beta_{0})))x_{k1}|}{\sqrt{2I_{n_{i}}(\beta_{0})(1,1)\log\log I_{n_{i}}(\beta_{0})(1,1)}}\geq\frac{1}{2}  \qquad a.s.
\end{eqnarray}

Now we define a $p\times1$ vector $\tilde{\beta}_{n_{i}}$ which includes only one non-zero component in the first coordinate, for $j=2,3,\cdots,p$
\begin{eqnarray}\label{3.43}
\tilde{\beta}_{n_{i}}(j)=\beta_{0j}, \quad {{\tilde \beta }_{{n_i}}}(1) = {C_0}\sqrt {\frac{{2\log \log {I_{{n_i}}}({\beta _0})(1,1)}}{{{I_{{n_i}}}({\beta _0})(1,1)}}}  + {\beta _{01}}
\end{eqnarray}
where ${C_0}$ is a constant which will be specified later.

Then, when $n_{i}$ is large enough,
\begin{eqnarray}\label{3.45}\nonumber
R_{2n_{i}}(\tilde{\beta}_{n_{i}})&=&\sum\limits_{k=1}^{n_{i}}w_{k}\dot{u}(x_{k}^{t}\beta_{0})(\dot{b}(x_{k}^{t}\beta_{0})-y_{k})x_{k}^{t}(\tilde{\beta}_{n_{i}}-\beta_{0})\nonumber\\
&=&\sum\limits_{k=1}^{n_{i}}w_{k}\dot{u}(x_{k}^{t}\beta_{0})(\dot{b}(x_{k}^{t}\beta_{0})-y_{k})x_{k1}(\tilde{\beta}_{n_{i}}(1)-\beta_{01})
\end{eqnarray}
Using the definition of $\tilde{\beta}_{n_{i}}$,  condition (H.1)(ii) and the \eqref{3.42}, we have
\begin{eqnarray}\label{3.46}\nonumber
R_{2n_{i}}(\tilde{\beta}_{n_{i}})&\leq& - \frac{1}{2}\sqrt {2{I_{{n_i}}}({\beta _0})(1,1)\log \log {I_{{n_i}}}({\beta _0})(1,1)}  \cdot {C_0}\sqrt {\frac{{2\log \log {I_{{n_i}}}({\beta _0})(1,1)}}{{{I_{{n_i}}}({\beta _0})(1,1)}}} \\\nonumber
&=& - {C_0}\log \log {I_{{n_i}}}({\beta _0})(1,1)\\
&\leq& - {C_0}{C_1}\log \log n_{i} \qquad a.s..
\end{eqnarray}
Where $C_1$ is the determined constant.
For ${R_{3{n_i}}}({{\tilde \beta }_{{n_i}}})$, we can see that
\begin{eqnarray}\nonumber
{R_{3{n_i}}}({{\tilde \beta }_{{n_i}}}) &= \sum\limits_{k = 1}^{{n_i}} {{w_k}\ddot u(x_k^t{\beta _1})} (\dot b(x_k^t{\beta _0}) - {y_k}){[x_{k1}^{}({{\tilde \beta }_{{n_i}}}(1) - {\beta _{01}})]^2}\\
& = C_0^2\frac{{2\log \log {I_{{n_i}}}({\beta _0})(1,1)}}{{{I_{{n_i}}}({\beta _0})(1,1)}}\sum\limits_{k = 1}^{{n_i}} {{w_k}\ddot u(x_k^t{\beta _1})} (\dot b(x_k^t{\beta _0}) - {y_k})x_{k1}^2
\end{eqnarray}
By using Corollary {\ref{cor 2}}, which implies
\begin{eqnarray}\label{LIL3}
\mathop {\lim \sup }\limits_{{n_i} \to \infty } \frac{{|\sum\limits_{k = 1}^{{n_i}} {{w_k}} \ddot u(x_k^t{\beta _1})[{y_k} - \dot b(u(x_k^t{\beta _0}))]x_{k1}^2|}}{{\sqrt {2{n_i}\log \log {n_i}} }} <\infty.
\end{eqnarray}

We also have sub-sequence $\{ {{n'}_i}{\rm{\} }}$ of $\{ {{n}_i}{\rm{\} }}$ such that
\begin{eqnarray}\label{3.42B}
\frac{{|\sum\limits_{k = 1}^{{{n'}_i}} {{w_k}} \ddot u(x_k^t{\beta _1})[{y_k} - \dot b(u(x_k^t{\beta _0}))]x_{k1}^2|}}{{\sqrt {2{{n'}_i}\log \log {{n'}_i}} }} \ge O(\frac{1}{{\rm{2}}}). \qquad a.s..
\end{eqnarray}

Analogous to the treatment of $R_{2n_{i}}(\tilde{\beta}_{n_{i}})$, we have
\begin{eqnarray}\label{3.46}\nonumber
R_{3n_{i}}(\tilde{\beta}_{n_{i}})&\leq&  - O(\frac{1}{2})\sqrt {2{{n'}_i}\log \log {{n'}_i}}  \cdot C_0^2\frac{{2\log \log {{n'}_i}}}{{{{n'}_i}}}  \nonumber\\
&=& - O(1)C_{0}\log \log {{n'}_i} \cdot \sqrt {\frac{{\log \log {{n'}_i}}}{{{{n'}_i}}}}  \qquad a.s..
\end{eqnarray}

Then we find that the sum of  ${R_{2{{n'}_i}}}({{\tilde \beta }_{{{n'}_i}}})$ and ${R_{3{{n'}_i}}}({{\tilde \beta }_{{{n'}_i}}})$ is
\begin{eqnarray}\label{sum}
{R_{2{n_i}}}({{\tilde \beta }_{{{n'}_i}}}) + {R_{3{n_i}}}({{\tilde \beta }_{{{n'}_i}}}) \le  - {C_0}(O(\sqrt {\frac{{\log \log {{n'}_i}}}{{{{n'}_i}}}} ) + C_{1})\log \log {{n'}_i}.
\end{eqnarray}

According to the $O(n)$ in \eqref{3.36}, set $O({{n'}_i})\leq C_2 {{n'}_i}$, and $C_{2}$ is a constant, we have
\begin{eqnarray}\label{r1}\nonumber
R_{1}(\tilde{\beta}_{{n'}_{i}})&\leq&{C_2}{n'}_{i}{\| {{{\tilde \beta }_{{{n'}_i}}} - {\beta _0}} \|^2}\\\nonumber
&=& O(1)\frac{{2{C_2}{C_0^2}{n'}_{i}}}{{n'}_{i}}\log \log {n'}_{i} \\
&\leq&C_{3}C_0^2\log \log {n'}_{i}.
\end{eqnarray}
Where $C_{3}$ is a constant. Therefore, according to (\ref{sum}) and (\ref{r1}), when ${n'}_{i}$ is sufficiently large, we have
\begin{eqnarray}\label{3.52}\nonumber
{K_{{n'}_{i}}}({{\tilde \beta }_{{n_i}}},{\beta _0})&=&R_{1}(\tilde{\beta}_{{n'}_{i}})+R_{2}(\tilde{\beta}_{{n'}_{i}})+R_{3}(\tilde{\beta}_{{n'}_{i}})\nonumber\\
&\leq& [{C_3}C_0^2 - {C_0}(O(\sqrt {\frac{{\log \log {{n'}_i}}}{{{{n'}_i}}}} ) + C_{1})]\log \log {{n'}_i}.
\end{eqnarray}
Setting ${C_3}C_0^2 - {C_0}(O(\sqrt {\frac{{\log \log {{n'}_i}}}{{{{n'}_i}}}} ) + C_1) < 0$, then ${C_0} < \frac{{O(\sqrt {\frac{{\log \log {{n'}_i}}}{{{{n'}_i}}}} ) + C_1}}{{{C_3}}}$. So we choose ${C_0}$ to be a certain positive constant which is smaller than $\frac{{O(\sqrt {\frac{{\log \log {{n'}_i}}}{{{{n'}_i}}}} ) + C_1}}{{{C_3}}}$. Thus,
\begin{eqnarray}\label{3.53}
{K_{{n'}_{i}}}({{\tilde \beta }_{{{n'}_{i}}}},{\beta _0}) \le  - O(1)\log \log {{n'}_{i}} \qquad a.s.
\end{eqnarray}
Notice the definition of $\hat{\beta}$ ($\hat \beta$ minimizes ${K_n}(\hat \beta )$), which ensures that ${K_{{n'}_{i}}}(\hat \beta ,{\beta _0}) \le {K_{{n'}_{i}}}({{\tilde \beta }_{{{n'}_{i}}}},{\beta _0})$. However when $n_{i}$ is sufficiently large, by assumption \eqref{3.35}, we get
\[{K_{{n'}_{i}}}({{\tilde \beta }_{{{n'}_{i}}}},{\beta _0})\le  - O(1)\log \log {{n'}_{i}} < o(1)\log \log {{n'}_{i}} = {K_{{n'}_{i}}}(\hat \beta ,{\beta _0})\]
via \eqref{smallo}. This is a contradiction, thus the proof of Theorem {\ref{theorem 1}} is completed. Therefore, for any $\alpha$ sub-model in $\Gamma_{c}$, $\mathop {\lim \sup }\limits_{n \to \infty } \frac{{\| {\hat \beta (\alpha ) - {\beta _0}(\alpha )}\|}}{{\sqrt {{n^{ - 1}}\log \log n} }}=d \ne 0$ \quad a.s. is true.

\section{Proof of Theorem {\ref{theorem 2}}}

According to the proof of Theorem {\ref{theorem 1}}, from the definition of ${K_n}(\beta ,{\beta _0})$, in order to prove (\ref{3.4}) for any $\alpha$ sub-model in $\Gamma_{c}$, it is needed to prove the following
\begin{eqnarray}\label{3.55}
0 \le {K_n}(\beta (\alpha ),{{\hat \beta }_0}(\alpha )) &=& {R_{1n}}(\hat \beta (\alpha )) + {R_{2n}}(\hat \beta (\alpha )) + {R_{3n}}(\hat \beta (\alpha ))\\
&=& O(\log \log n)~a.s..
\end{eqnarray}
Because $\hat{\beta}$ is the MLE of $\beta$, then ${l_n}(\beta_{0}(\alpha)) \le {l_n}(\hat \beta (\alpha ))$, that is,
\begin{eqnarray}\nonumber
{K_n}(\beta (\alpha ),{{\hat \beta }_0}(\alpha )) = {l_n}(\hat \beta (\alpha )) - {l_n}(\beta ({\alpha _0})) \ge 0 \qquad a.s..
\end{eqnarray}
Notice Lemma {\ref{lemma a.1}} and Theorem {\ref{theorem 1}}
\begin{eqnarray}\label{R1}\nonumber
|{R_{1n}}(\hat \beta (\alpha ))|&=&|\sum\limits_{k=1}^{n}w_{k}K(u(x_{k\alpha}^{t}\beta(\alpha )),u(x_{k\alpha}^{t}\beta_{0}(\alpha )))|\nonumber\\
&\leq& |\frac{1}{2}C_{u}\sum\limits_{k=1}^{n}w_{k}{{\dot u}^2}(x_{k\alpha}^t{\beta_{1}(\alpha ) })|(x_{k\alpha}^{t}\beta(\alpha)-x_{k\alpha}^{t}\beta_{0}(\alpha))^{2}\nonumber\\
&=&O(n)\|\beta(\alpha )-\beta_{0}(\alpha)\|^{2}\nonumber\\
&=&O(\log\log n).
\end{eqnarray}
Then, according to Corollary {\ref{cor 1}} and Theorem {\ref{theorem 1}}, we have
\begin{eqnarray}\label{R2}\nonumber
|R_{2n}(\hat{\beta}(\alpha ))|&=&|\sum\limits_{k = 1}^n {{w_k}} \dot u(x_{k\alpha }^t\beta_{0}(\alpha))[{y_k} - \dot b(u(x_{k\alpha }^t\beta_{0}(\alpha)))]{x_{k\alpha }}|\cdot\|\hat{\beta}(\alpha )-\beta_{0}(\alpha)\|\nonumber\\
&=&O(\sqrt{n\log\log n})\cdot O(\sqrt{n^{-1}\log\log n})=O(\log\log n) ~a.s..
\end{eqnarray}
For ${R_{3n}}(\hat \beta (\alpha ))$, by LIL in Corollary {\ref{cor 2}}, we get almost surely
\begin{eqnarray}\label{LIL4}
\left| {\sum\limits_{k = 1}^n {{{w}_{k}}} \ddot u(x_{k\alpha }^t\beta ({\alpha _0}))[{y_k} - \dot b(u(x_{k\alpha }^t\beta ({\alpha _0})))]||{x_{k\alpha }}||^2} \right| = O(\sqrt {n\log \log n} ).
\end{eqnarray}
Therefore
\begin{eqnarray}\label{R3}\nonumber
|{R_{3n}}(\hat \beta (\alpha ))|&=&|\sum\limits_{k = 1}^n {{w}_{k}} \dot u(x_{k\alpha }^t\beta_{0}(\alpha))[{y_k} - \ddot b(u(x_{k\alpha }^t\beta_{0}(\alpha)))]||{x_{k\alpha }}||^2|\cdot\|\hat{\beta}(\alpha )-\beta_{0}(\alpha)\|_{\infty}\nonumber\\
&=&  O(\sqrt{n\log\log n})\cdot O(\sqrt{n^{-1}\log\log n})=O(\log\log n) \qquad a.s..
\end{eqnarray}

Combining the estimates (\ref{R1}), (\ref{R2}) and (\ref{R3}), we obtain
\begin{eqnarray}\nonumber
|{K_n}(\beta (\alpha ),{{\hat \beta }_0}(\alpha ))| \le |{R_{1n}}(\hat \beta (\alpha ))| + |{R_{2n}}(\hat \beta (\alpha ))| + |{R_{3n}}(\hat \beta (\alpha ))| = O(\log \log n)\quad a.s..
\end{eqnarray}
Therefore, Theorem {\ref{theorem 2}} is proved.

\section{Proof of Theorem {\ref{theorem 3}}}

First, let ${{\hat \beta }^*}(\alpha )$ be the $p_{\alpha}$-dimensional vector which is defined by augmenting ${{\hat \beta }}(\alpha )$ with  $p-p_{\alpha}$ 0's such that the sub-vector of ${{\hat \beta }^*}(\alpha )$ indexed by $\alpha$ matches ${{\hat \beta }}(\alpha )$. Then, it is easy to see that proving the (\ref{3.5}) is tantamount to give the proof of the following argument: for any incorrect model $\alpha\in \Gamma_{w}$
\begin{eqnarray}\label{3.60}
\liminf\limits_{n\rightarrow\infty}n^{-1}{K_n}({{\hat \beta }^*}(\alpha ),{\beta _0})>0 \qquad a.s..
\end{eqnarray}

Next, we define $l_2$ ball with radius $a:=\frac{1}{2}\min\limits_{1\leq i\leq p_{\alpha_{0}}}|\beta_{0}(\alpha_{0})_{i}|$,
\begin{eqnarray}\nonumber
B_{0}=\{\beta:\|\beta-\beta_{0}\|< a\},
\end{eqnarray}
where $\alpha_{0}$ is the minimum dimension of the true model belonging to $\Gamma_{c}$, and $\beta_{0}(\alpha_{0})_{i}$ is the $i$th component of ${K_n}(\hat \beta ,{\beta _0})$.

Obviously, $B_{0}$ is a tight set. By the definition of incorrect model $\alpha\in \Gamma_{w}$, we get $\|\hat{\beta}^{*}(\alpha)-\beta_{0}\|\geq \min\limits_{1\leq i\leq p_{\alpha_{0}}}|\beta_{0}(\alpha_{0})_{i}|$. Then we have $\hat{\beta}^{*}(\alpha)\notin B_{0}$. By applying Theorem {\ref{theorem 1}}, for the large $n$, MLE $\hat{\beta}$ is almost surely an interior point of $B_{0}$. Applying the convexity of ${K_n}(\beta ,{\beta _0})$, we have
\begin{eqnarray}\label{3.61}
{K_n}({{\hat \beta }^*}(\alpha ),{\beta _0}) \ge \mathop {\inf }\limits_{\beta  \in \partial {B_0}} {K_n}(\beta ,{\beta _0}),
\end{eqnarray}
where $\partial B_{0}=\{\beta:\|\beta-\beta_{0}\|=\frac{1}{2}\min\limits_{1\leq i\leq p_{\alpha_{0}}}|\beta_{0}(\alpha_{0})_{i}|=a\}$,  which is the boundary of $B_{0}$.

To prove (\ref{3.60}), it is required to prove the following result:
\begin{eqnarray}\label{3.62}
\liminf\limits_{n\rightarrow\infty}\inf\limits_{\beta\in\partial B_{0}}n^{-1}{K_n}(\beta ,{\beta _0})>0 \qquad a.s..
\end{eqnarray}
Using the similar argument in \eqref{3.26b}, we have
\begin{eqnarray}\label{3.63}\nonumber
{R_{1n}}(\beta)I(\beta\in\partial B_{0})
&\ge &O(n)\|\beta-\beta_{0}\|^{2}I(\beta\in\partial B_{0}).
\end{eqnarray}
Then
\begin{eqnarray}\label{3.64}
\mathop {\inf }\limits_{\beta  \in \partial {B_0}} {R_{1n}}(\beta) \ge \mathop {\inf }\limits_{\beta  \in \partial {B_0}} O(n){\left\| {\beta  - {\beta _0}} \right\|^2} = O(n)a^{2} = O(n).
\end{eqnarray}
Similarly, using (\ref{3.28}) one has
\begin{eqnarray}\label{3.65}
|{R_{2n}}(\beta)|I(\beta\in\partial B_{0})&\leq&\|\sum\limits_{k=1}^{n}w_{k}\dot{u}(x_{k}^{t}\beta_{0})(y_{k}-\dot{b}(u(x_{k}^{t}\beta_{0})))x_{k}\|\cdot\|\beta-\beta_{0}\|I(\beta\in\partial B_{0})\nonumber\\
&=&O(\sqrt{n\log\log n})\|\beta-\beta_{0}\|  \qquad a.s..
\end{eqnarray}
Then by taking supreme, we get
\begin{eqnarray}\label{3.66}
\sup\limits_{\beta\in\partial B_{0}}|{R_{2n}}(\beta)|=O(\sqrt{n\log\log n}) \qquad a.s..
\end{eqnarray}
For, using $\eqref{R3S}$ we similarly have
\begin{eqnarray}\label{3.67}\nonumber
|{R_{3n}}(\beta)|I(\beta\in\partial B_{0})\leq \| {\beta  - {\beta _0}}\| O(\sqrt {n\log \log n} )\qquad a.s..
\end{eqnarray}
By taking supreme, it gives
\begin{eqnarray}\label{3.69}
\sup\limits_{\beta\in\partial B_{0}}|{R_{3n}}(\beta)|\leq a O(\sqrt {n\log \log n} )=O(\sqrt {n\log \log n} ) \qquad a.s..
\end{eqnarray}
According to (\ref{3.64}), \eqref{3.66} and \eqref{3.69}, we obtain
\begin{eqnarray}
\inf\limits_{\beta\in\partial B_{0}}n^{-1}{K_n}(\beta ,{\beta _0})\geq\inf\limits_{\beta\in\partial B_{0}}R_{1n}(\beta,n)-\sup\limits_{\beta\in\partial B_{0}}|R_{2n}(\beta,n)|-\sup\limits_{\beta\in\partial B_{0}}|R_{3n}(\beta,n)|.\nonumber\\
~~\label{3.70}
\end{eqnarray}
Then
\begin{eqnarray}\nonumber
\inf\limits_{\beta\in\partial B_{0}}n^{-1} {K_n}(\beta_0 ,{\beta})&\geq&O(1)-O(\sqrt{n^{-1}\log\log n}).
\end{eqnarray}
Therefore, when $n\rightarrow\infty$, $\inf\inf\limits_{\beta\in\partial B_{0}}n^{-1}{K_n}(\beta_0 ,{\beta })>0$. Hence (\ref{3.62}) and (\ref{3.60}) is established. Then Theorem {\ref{theorem 3}} is proved.

\section{Proof of Theorem {\ref{theorem 4}}}

The proof is to divide $\alpha$ into two cases. That is: 1.$\alpha\in \Gamma_{c}$; 2.$\alpha\in \Gamma_{w}$.

Step1: For any correct model $\alpha\in \Gamma_{c}$, employing Theorem {\ref{theorem 2}} we have
\begin{eqnarray}\nonumber
{S_n}(\alpha ): &=& - {l_n}(\hat \beta (\alpha )) + C(n,\hat \beta (\alpha ))\nonumber\\
&=&- {l_n}(\beta ({\alpha _0})) + C(n,\hat \beta (\alpha )) - O(\log \log n) \qquad  a.s.,
\end{eqnarray}
Minimizing over $\alpha$ gives
\[\mathop {\min }\limits_{\alpha  \in {\cal A}} {S_n}(\alpha ) =  - {l_n}(\beta ({\alpha _0})) + \mathop {\min }\limits_{\alpha  \in {\cal {\cal A}}} [C(n,\hat \beta (\alpha )) - O(\log \log n)] \qquad  a.s.,\]
Since the right hand side above should not be a decreasing function as $n$ is large, it leads to $C(n,\hat \beta (\alpha )) > O(\log \log n)$. The reason is that for any correct model $\alpha\in \Gamma_{c}$ we need following equality:
\[\mathop {\min }\limits_{\alpha  \in {\cal A}} [C(n,\hat \beta (\alpha )) - O(\log \log n)] = \mathop {\min }\limits_{\alpha  \in {\cal A}} C(n,\hat \beta (\alpha ))\qquad  a.s.\]
Note that the definition of  ${\alpha _0}$ is ${\alpha _0} := \mathop {\arg \min}\limits_{\alpha  \in A,\alpha  \in {\Gamma _c}} {p_\alpha }$. Then, by the definition of ${\hat \alpha }$ and the fact that $C(n,\hat \beta (\alpha ))$ is an increasing function of ${p_\alpha }$, we have
\[\hat \alpha  = \mathop {\arg \min }\limits_{\alpha  \in {\cal A}} {S_n}(\alpha ) = \mathop {\arg \min }\limits_{\alpha  \in {\cal A}} C(n,\hat \beta (\alpha )) = \mathop {\arg \min }\limits_{\alpha  \in {\cal A},\alpha  \in {\Gamma _c}} {p_\alpha }\qquad  a.s..\]
Thus $\hat \alpha  = {\alpha _0}$~a.s..

Step2: For any wrong model $\alpha\in \Gamma_{w}$, employing Theorem {\ref{theorem 3}} we get
\[{S_n}(\alpha ): =  - {l_n}(\beta ({\alpha _0})) + C(n,\hat \beta (\alpha )) + O(n).\]
Since $\alpha$ is wrong, thus ${\alpha _0} \ne \hat \alpha$, and then $\alpha _0$ should not be obtained by minimizing ${S_n}(\alpha )$. Therefore ${S_n}(\alpha )$ will increase as $\alpha$ approximates $\alpha _0$. Observe that minimizing ${S_n}(\alpha )$ is equivalent to minimizing $C(n,\hat \beta (\alpha )) + O(n)$. Consequently, if we assume that minimizing $C(n,\hat \beta (\alpha )) + O(n)$ is the same as minimizing $C(n,\hat \beta (\alpha ))$, then this assumption leads to $\hat \alpha = {\alpha _0}$ a contradiction with ${\alpha _0} \ne \hat \alpha$. Hence, we must have $C(n,\hat \beta (\alpha )) < O(n)$.

In summary, if $O(\log \log n) < C(n,\hat \beta (\alpha )) < O(n)$, then the model selection criteria are strongly consistent. Finally, via the order of AIC penalty term, it is not strongly consistent. Similarly, both the order of BIC and SCC are $\log n$ and it follows that they are strongly consistent.

\section{Proof of Theorem {\ref{theorem 5}}}
Before the proof, we pose two LIL results for $\rho $-mixing process and $m$-dependent responses.
\begin{lemma}{(LIL of $\rho $-mixing process, Corollary~9.2.1 in \cite{Lin1997})}{\label{lemma h.1}}
 Let $\{Z_{n},n\geq1\}$ be a strictly
stationary $\rho$-mixing sequence with $E(Z_{1})=0$, $E(Z_{1}^{2})<\infty$. Let $S_{n}=\sum\limits_{k=1}^{n}Z_{k}$. Assume that\\
(i) $\sigma _n^2: = E(S_n^2) \to \infty $~\text{as}~$(n\rightarrow\infty)$;\\
(ii) $\rho (n) = O({(\log n)^{ - 1 - \varepsilon }})$~\text{for~some}~$\varepsilon  > 0$.\\
Then, we have
\begin{eqnarray}\nonumber
\limsup\limits_{n\rightarrow\infty}\frac{|S_{n}|}{\sqrt{2\sigma_{n}^{2}\log\log\sigma_{n}^{2}}}=1 \qquad a.s..
\end{eqnarray}
\end{lemma}

\begin{lemma}{(LIL of $m$-dependent random variables, \cite{Chen97})}\label{lemma h.2}
 Let $\{Z_{n},n\geq1\}$ be a real stationary strongly mixing sequences with $E(Z_{1})=0$, $E(Z_{1}^{2})<\infty$. If $\{Z_{n},n\geq1\}$ is $m$-dependent and let $S_{n}=\sum\limits_{k=1}^{n}Z_{k}$, we have
\begin{eqnarray}\nonumber
\limsup\limits_{n\rightarrow\infty}\frac{|S_{n}|}{\sqrt{2n \log\log n}}=\sigma^{2} \qquad a.s.,
\end{eqnarray}
where $\sigma _{}^2: = EZ_1^2 + 2\sum\limits_{k = 2}^{m + 1} {E{X_1}{X_k}}$.
\end{lemma}

The proof of \eqref{3.2wd} for $\rho $-mixing responses or $m$-dependent responses is slightly different from what we mentioned in Theorem {\ref{theorem 1}} concerning independence case, and the main dissimilarity is the calculating of the variance of the score function. We cannot ignore the non-zero cross terms which will be included as the covariance parts.

\textbf{(i): $\rho$-mixing sequence}:
For $j=1,\cdots,p$, let
\[
{S_n}{(\beta_0 )_k}: = \sum\limits_{k = 1}^n {{w_k}} {x_{kj}}\dot u(x_k^t\beta_0 )[{y_k} - \dot b(u(x_k^t\beta_0 ))] = \sum\limits_{k = 1}^n {{{\mathord{\buildrel{\lower3pt\hbox{$\scriptscriptstyle\smile$}}
\over w} }_k}} [{y_k} - E({y_k})].
\]
where ${{\mathord{\buildrel{\lower3pt\hbox{$\scriptscriptstyle\smile$}}
\over w} }_k} := {w_k}{x_{kj}}\dot u(x_k^t\beta _0)$.

We have
\begin{eqnarray}\nonumber
{I_n}(\beta )(k,k) = Var({S_n}{(\beta )_k}) &=& \sum\limits_{k = 1}^n {\mathord{\buildrel{\lower3pt\hbox{$\scriptscriptstyle\smile$}}
\over w} _k^2} Var({y_k}) +2 \sum\limits_{1 \le i < j \le n} {{{\mathord{\buildrel{\lower3pt\hbox{$\scriptscriptstyle\smile$}}
\over w} }_i}{{\mathord{\buildrel{\lower3pt\hbox{$\scriptscriptstyle\smile$}}
\over w} }_j}{\mathop{\rm Cov}} ({y_i},{y_j})}\nonumber \\
&=& O(n) + 2 \sum\limits_{1 \le i < j \le n} {{{\mathord{\buildrel{\lower3pt\hbox{$\scriptscriptstyle\smile$}}
\over w} }_i}{{\mathord{\buildrel{\lower3pt\hbox{$\scriptscriptstyle\smile$}}
\over w} }_j}{\mathop{\rm Cov}} ({y_i},{y_j})} \label{IF}
\end{eqnarray}

\begin{lemma}{(Davydov's inequality)}\label{lemma h.3}
 Let $X \in \mathcal{F}_1^k, Y \in \mathcal{F}_{k + n}^\infty$ with $E|X|^{p}<\infty$ and $E|Y|^{p}<\infty$~(${p^{ - 1}} + {q^{ - 1}} < 1$), then
\[\left| {EXY - EXEY} \right| \le 10\sqrt[p]{{E{|X|^p}}}\sqrt[q]{{E{|Y|^q}}}{(\rho (n))^{1 - 1/p - 1/q}},\]
(see Lemma~1.2.4 in \cite{Lin1997} or Corollary~1.1.1 in \cite{Bosq1998}).
\end{lemma}
By Davydov's inequality  above, let $p = q = 3$ in Lemma {\ref{lemma h.3}}, we have
$$| {\sum\limits_{1 \le i < j \le n} {{{\mathord{\buildrel{\lower3pt\hbox{$\scriptscriptstyle\smile$}}
\over w} }_i}{{\mathord{\buildrel{\lower3pt\hbox{$\scriptscriptstyle\smile$}}
\over w} }_j}{\mathop{\rm cov}} ({y_i},{y_j})} } | \le \sum\limits_{m = 1}^{n - 1} {\sum\limits_{1 \le i < j \le n,i - j = m} {10{{(\alpha (n))}^{1/3}}(E|y_i|^3E|y_j|^3})^{1/3}|{\mathord{\buildrel{\lower3pt\hbox{$\scriptscriptstyle\smile$}}
\over w} _i}{\mathord{\buildrel{\lower3pt\hbox{$\scriptscriptstyle\smile$}}
\over w} _j}| } .$$
By (H.5): The error sequence ${\varepsilon _i}$ satisfies $\rho$-mixing condition with geometric decay: $\rho (m) =O({r^{ - m}})$, and the relationship between $\alpha$-mixing coefficient and $\rho $-mixing coefficient meets: $\alpha (n) \le \frac{{\rm{1}}}{{\rm{4}}}\rho (n)$, we have
\begin{eqnarray}\nonumber
&~&| {\sum\limits_{1 \le i < j \le n} {{{\mathord{\buildrel{\lower3pt\hbox{$\scriptscriptstyle\smile$}}
\over w} }_i}{{\mathord{\buildrel{\lower3pt\hbox{$\scriptscriptstyle\smile$}}
\over w} }_j}{\mathop{\rm cov}} ({y_i},{y_j})} } |\nonumber\\
&=& O(1)\sum\limits_{m = 1}^{n - 1} {\sum\limits_{1 \le i < j \le n,i - j = m} {\rho{{(m)}^{1/3}}} }  = O(1)\sum\limits_{m = 1}^{n - 1} {\sum\limits_{1 \le i < j \le n,i - j = m} {O({{({r^{1/3}})}^{ - n}})} } \nonumber\\
&=& O(1)[(n - 1){({r^{1/3}})^{ - 1}} + (n - 2){({r^{1/3}})^{ - 2}} +  \cdots  + (n - (n - 1)){({r^{1/3}})^{ - (n - 1)}}]\nonumber\\
&=& O(1)[n\sum\limits_{m = 1}^{n - 1} {{{({r^{1/3}})}^{ - m}}}  - \sum\limits_{m = 1}^{n - 1} {m{{({r^{1/3}})}^{ - m}}} ]\nonumber\\
&=& O(1)[nO(1) - O(1)] = O(n).\label{rho}
\end{eqnarray}

Then \eqref{IF} and \eqref{rho} imply ${I_n}(\beta )(k,k) = O(n) \to \infty $. So condition (i) in Lemma {\ref{lemma h.1}} is satisfied. And condition (ii) in Lemma {\ref{lemma h.1}} is also valid since the (H.5) gives $\rho (n) = O({r^{ - n}}) = O({(\log n)^{ - 1 - \varepsilon }})$.

\textbf{(ii) $m$-dependent sequence}:
For $j=1,\cdots,p$, we have
\begin{align*}
{I_n}(\beta )(k,k) &= O(n) + 2\sum\limits_{1 \le i < j \le n} {{{\mathord{\buildrel{\lower3pt\hbox{$\scriptscriptstyle\smile$}}
\over w} }_i}{{\mathord{\buildrel{\lower3pt\hbox{$\scriptscriptstyle\smile$}}
\over w} }_j}{\mathop{\rm cov}} ({y_i},{y_j})}\\
=& O(n) + O(1)\sum\limits_{z = 1}^m {\sum\limits_{1 \le i < j \le n,i - j = z} {{{\mathord{\buildrel{\lower3pt\hbox{$\scriptscriptstyle\smile$}}
\over w} }_i}{{\mathord{\buildrel{\lower3pt\hbox{$\scriptscriptstyle\smile$}}
\over w} }_j}{\mathop{\rm cov}} ({y_i},{y_j})} }  = O(n).
\end{align*}
Under (H.6), then we can use a result of the LIL for real stationary strongly mixing sequences with $m$-dependent properties, see Lemma {\ref{lemma h.2}}.

In summary, the LIL for the $\alpha$-mixing or $m$-dependent score functions is
\[\mathop {\lim \sup }\limits_{n \to \infty } \frac{{|\sum\limits_{k = 1}^n {{w_k}} {x_{kj}}\dot u(x_k^t\beta _0)[{y_k} - \dot b(u(x_k^t\beta_0 ))]|}}{{\sqrt {2n\log \log n} }} = O(1) \qquad  a.s..\]
Therefore, the proof of Theorem {\ref{theorem 5}} can be imitated from the proof of Theorem {\ref{theorem 1}}.

After Observing the proofs of Theorem {\ref{theorem 2}}, Theorem {\ref{theorem 3}} and Theorem {\ref{theorem 4}} do not involve the calculation of the second moment of the weighted score function, we conclude that the corresponding new proofs are consistent with the proofs of the case of independent responses.

We only pay attention to the equations which concern the LIL of the several weighted score functions, i.e. \eqref{LIL1}, \eqref{LIL2} and \eqref{LIL3}. Next, for weakly dependent version of Theorem {\ref{theorem 3}}, the LIL \label{LIL4} is also true when the weighted score function is weakly dependent sums.  As for Theorem {\ref{theorem 3}}, the proof is the same by applying the LIL of weakly dependent weighted score function.

Then as for two dependent cases, assuming that conditions (H.1) to (H.6) are satisfied, we also have the same conclusions as in Theorem {\ref{theorem 2}} and Theorem {\ref{theorem 3}}. The proofs in other places in Theorem {\ref{theorem 6}} are consistent with Theorem {\ref{theorem 4}}. Therefore we do not make repetitions here.
\end{appendix}


\begin{thebibliography}{999}

\bibitem[Ai et al.(2020)]{Ai2018}
Ai, M., Yu, J., Zhang, H., \& Wang, H. (2020). Optimal Subsampling Algorithms for Big Data Regressions.  Statistica Sinica. DOI: 10.5705/ss.202018.0439

\bibitem[{Akaike(1973)}]{Akaike1973}
Akaike, H. (1973). Information theory and an extension of the maximum likelihood principle. In 2nd ed International Symposium on Information Theory, 1973. Akademiai Kaido.

\bibitem[{Bosq(1998)}]{Bosq1998}
Bosq, D. (1998). Nonparametric statistics for stochastic processes: estimation and prediction. Springer.

\bibitem[{Brown(1986)}]{Brown1986}
Brown, L. D. (1986). Fundamentals of statistical exponential families: with applications in statistical decision theory. IMS.

\bibitem[{Chen(1997)}]{Chen97}
Chen, X. (1997). The law of the iterated logarithm for m-dependent Banach space valued random variables. Journal of Theoretical Probability, 10(3), 695-732.

\bibitem[{Chen(2011)}]{Chen11}
Chen, X. (2011). Quasi Likelihood Method for Generalized Linear Model~(in Chinese). Press of University of Science and Technology of China.

\bibitem[{Czado and Munk (2000)}]{Czado00}
Czado, C., \& Munk, A. (2000). Noncanonical links in generalized linear models when is the effort justified?. Journal of Statistical Planning and Inference, 87(2), 317-345.

\bibitem[Efron and Hastie(2016)]{Efron2016}
Efron, B., \& Hastie, T. (2016). Computer age statistical inference: algorithms, evidence, and data science. Cambridge University Press.

\bibitem[{Fahrmeir and Kaufmann(1985)}]{Fahrmeir85}
Fahrmeir, L., \& Kaufmann, H. (1985). Consistency and asymptotic normality of the maximum likelihood estimator in generalized linear models. The Annals of Statistics, 342-368.

\bibitem[{Fahrmeir and Tutz(2001)}]{Fahrmeir2001}
Fahrmeir, L., \& Tutz, G. (2001). Multivariate statistical modelling based on generalized linear models, 2ed. Springer.

\bibitem[{Fan et al.(2016)}]{Fan16}
Fan, J., Qi, L., \& Tong, X. (2016). Penalized least squares estimation with weakly dependent data. Science China Mathematics, 59(12), 2335-2354.

\bibitem[{Fang(1998)}]{Fang1998}
Fang, X. (1998). Laws of the iterated logarithm for maximum likelihood estimates of parameter vectors in nonhomogeneous Poisson processes (Chinese). Acta Scientiarum Naturalium Universitatis Pekinensis, 1998, 34(5), 563-573. MR1681923

\bibitem[{Hansen(2018)}]{Hansen18}
Hansen, B. (2018). Econometrics. Version: Jan 2018\\ \url{https://www.ssc.wisc.edu/~bhansen/econometrics}.

\bibitem[He and Wang(1995)]{He1995}
He, X., \& Wang, G. (1995). Law of the iterated logarithm and invariance principle for M-estimators. Proceedings of the American Mathematical Society, 123(2), 563-573.

\bibitem[Kim and Jeon(2016)]{Kim2016}
Kim, Y., \& Jeon, J. J. (2016). Consistent model selection criteria for quadratically supported risks. The Annals of Statistics, 44(6), 2467-2496.

\bibitem[{Kroll(2019)}]{Kroll2018}
Kroll, M. (2019). Non-parametric Poisson regression from independent and weakly dependent observations by model selection. Journal of Statistical Planning and Inference, 199, 249-270.

\bibitem[Lai and Wei(1982)]{Lai1982}
Lai, T. L., \& Wei, C. Z. (1982). A law of the iterated logarithm for double arrays of independent random variables with applications to regression and time series models. The Annals of Probability, 320-335.

\bibitem[Lin and Lu(1997)]{Lin1997}
Lin, Z., \& Lu, C. (1997). Limit theory for mixing dependent random variables. Springer.

\bibitem[{Markatou et al.(1998)}]{Markatou98}
Markatou, M., Basu, A., \& Lindsay, B. G. (1998). Weighted likelihood equations with bootstrap root search. Journal of the American Statistical Association, 93(442), 740-750.

\bibitem[{McCullagh and Nelder(1989)}]{McCullagh1989}
McCullagh, P. \& Nelder, J. A. (1989). Generalized linear models. Second Edition. Chapman and Hall, London.

\bibitem[Miao and Yang(2011)]{Miao2011}
Miao, Y., \& Yang, G. (2011). The loglog law for LS estimator in simple linear EV regression models. Statistics, 45(2), 155-162.

\bibitem[{Nelder(1972)}]{Nelder1972}
Nelder, J., \& Wedderburn, R. (1972). Generalized Linear Models. Journal of the Royal Statistical Society. Series A (General), 135(3), 370-384.

\bibitem[{Petrov(1995)}]{Petrov1995}
Petrov, V. V. (1995). Limit theorems of probability theory: sequences of independent random variables. Oxford, New York.


\bibitem[Qian and Wu(2006)]{Qian2006}
Qian, G., \& Wu, Y. (2006). Strong limit theorems on model selection in generalized linear regression with binomial responses. Statistica Sinica, 1335-1365.

\bibitem[Rao and Wu(1989)]{Rao1989}
Rao, R., \& Wu, Y. (1989). A strongly consistent procedure for model selection in a regression problem. Biometrika, 76(2), 369-374.

\bibitem[Rao and Zhao(1992)]{Rao1992}
Rao, C. R., \& Zhao, L. C. (1992). Linear representation of M-estimates in linear models. Canadian Journal of Statistics, 20(4), 359-368.

\bibitem[{Rigollet(2012)}]{Rigollet12}
Rigollet, P. (2012). Kullback-Leibler aggregation and misspecified generalized linear models. The Annals of Statistics, 40(2), 639-665.

\bibitem[Rissanen(1989)]{Rissanen1989}
Rissanen, J. (1989). Stochastic complexity in statistical inquiry. World Scientific.

\bibitem[{Schwarz(1978)}]{Schwarz1978}
Schwarz, G. (1978). Estimating the dimension of a model. The Annals of Statistics, 6(2), 461-464.

\bibitem[Stout(1974)]{Stout1974}
Stout, W. F. (1974). Almost sure convergence, Academic Press, NewYork.

\bibitem[Shao(2003)]{Shao2003}
Shao, J. (2003). Mathematical Statistics 2ed. Springer New York.

\bibitem[Tutz(2011)]{Tutz2011}
Tutz, G. (2011). Regression for categorical data (Vol. 34). Cambridge University Press.

\bibitem[{van der Vaart(1998)}]{Vaart1998}
van der Vaart, A. W. (1998). Asymptotic statistics (Vol. 3). Cambridge university press.

\bibitem[{Wright(2017)}]{Wright2017}
Wright, S. J. (2017). Optimization algorithms for data analysis. IAS/Park City Mathematics Series, to appear.

\bibitem[Wu and Zen(1999)]{Wu1999}
Wu, Y., \& Zen, M. M. (1999). A strongly consistent information criterion for linear model selection based on M-estimation. Probability Theory and Related Fields, 113(4), 599-625.

\bibitem[Yin et al.(2006)]{Yin2006}
Yin, C., Zhao, L., \& Wei, C. (2006). Asymptotic normality and strong consistency of maximum quasi-likelihood estimates in generalized linear models. Science in China Series A, 49(2), 145-157.

\bibitem[Zhang and Jia(2017)]{Zhang2017}
Zhang, H., \& Jia, J. (2017). Elastic-net regularized high-dimensional negative binomial regression: consistency and weak signals detection. arXiv preprint arXiv:1712.03412.

\bibitem[Zheng and Peng(2017)]{Zheng2017}
Zheng, Q., \& Peng, L. (2017). Consistent model identification of varying coefficient quantile regression with BIC tuning parameter selection. Communications in Statistics-Theory and Methods, 46(3), 1031-1049.
\end{thebibliography}
\end{document}